\documentclass[11pt]{article}
\usepackage{geometry,amsfonts,epsf,amsmath,amssymb,latexsym,graphicx}
\usepackage[utf8]{inputenc}
\usepackage{tikz}

\newtheorem{thm}{\bf Theorem}[section]
\newtheorem{cor}[thm]{\bf Corollary}
\newtheorem{lemma}[thm]{\bf Lemma}
\newtheorem{proposition}[thm]{\bf Proposition}

\newtheorem{observation}[thm]{\bf Observation}

\newcommand{\proof}{\noindent{\bf Proof.\ }}
\newcommand{\qed}{\hfill $\square$ \bigskip}
\newcommand{\smallqed}{{\tiny ($\Box$)}}

\newcommand{\D}{{Dominator }}
\newcommand{\St}{{Staller }}

\newenvironment{unnumbered}[1]{\trivlist
\item [\hskip \labelsep {\bf #1}]\ignorespaces\it}{\endtrivlist}

\newcommand{\remove}[1]{}
\newcommand{\dstartpass}{\gamma_{g}^{{\rm sp}}}
\newcommand{\sstartpass}{\gamma_{g}'^{{\rm sp}}}
\newcommand{\xbar}{{\overline{x}}}
\newcommand{\uvprop}{$\{u,v\}$-inclusion property}
\newcommand{\Buvprop}{$\mathbf{\{u,v\}}$\textbf{-inclusion property}}

\let\oldenumerate\enumerate
\renewcommand{\enumerate}{
  \oldenumerate
  \setlength{\itemsep}{0pt}
  \setlength{\parskip}{0pt}
  \setlength{\parsep}{0pt}
}

\newcommand{\dstart}{\gamma_g}
\newcommand{\sstart}{\gamma_g^\prime}
\newcommand{\w}{{\rm w}}
\newcommand{\1}{ \vspace{0.1cm} }
\newcommand{\2}{ \vspace{0.2cm} }
\newcommand{\modo}{{\rm mod \,}}

\begin{document}

\title{Cutting Lemma and Union Lemma for the \\ Domination Game}

\author{$^1$Paul Dorbec, \, $^2$Michael A. Henning, \, $^{3,4,5}$Sandi Klav\v zar, \, $^{5,6}$Ga\v{s}per Ko\v{s}mrlj \\ \\
  $^1$Univ. Bordeaux, LaBRI, UMR5800, F-33405 Talence \\
  CNRS, LaBRI, UMR5800, F-33405 Talence \\
  Email: {\tt paul.dorbec@u-bordeaux.fr}   \\
  \\
  $^2$Department of Pure and Applied Mathematics\\
  University of Johannesburg \\
  Auckland Park, 2006 South Africa \\
  Email:  {\tt mahenning@uj.ac.za} \\
  \\
  $^3$Faculty of Mathematics and Physics, University of Ljubljana, Slovenia\\
  $^4$Faculty of Natural Sciences and Mathematics, University of Maribor, Slovenia\\
  $^5$Institute of Mathematics, Physics and Mechanics, Ljubljana, Slovenia\\
  Email:  {\tt sandi.klavzar@fmf.uni-lj.si} \\
  \\
  $^6$Abelium R\&D, Ljubljana, Slovenia \\
  Email:  {\tt gasperk@abelium.eu}
}

\date{}
\maketitle

\begin{abstract}
Two new techniques are introduced into the theory of the domination game. The cutting lemma bounds the game domination number of a partially dominated graph with the game domination number of suitably modified partially dominated graph.  The union lemma bounds the S-game domination number of a disjoint union of paths using appropriate weighting functions. Using these tools a conjecture asserting that the so-called three legged spiders are game domination critical graphs is proved. An extended cutting lemma is also derived and all game domination critical trees on 18, 19, and 20 vertices are listed.
\end{abstract}


\medskip\noindent
{\bf AMS Subj.\ Class:} 05C57, 05C69, 91A43

\section{Introduction}

The domination game was introduced in 2010~\cite{bresar-2010}. Until now, the game received considerable attention. In the last few years more than thirty papers have been published on the game and its total variation, the later being introduced in~\cite{henning-2015}. Central themes studied on these games are the $3/5$-Conjecture for the usual game~\cite{kinnersley-2013} and the $3/4$-Conjecture for the total game~\cite{henning-2017a}. These conjectures were investigated in~\cite{bresar-2013, bujtas-2015a, bujtas-2015b, henning-2016, henning-2017b, marcus-2016+, schmidt-2016} and~\cite{bujtas-2017+, henning-2017a, henning-2017+}, respectively. For additional aspects of the domination game we refer to recent papers~\cite{bresar-2016, bresar-2017, nadjafi-arani-2016, XLK-2018}.

Let $G$ be a graph. The domination game is played on $G$ by Dominator and Staller, who take turns choosing a vertex from $G$. Each vertex chosen must dominate at least one vertex not dominated by the set of vertices previously chosen. The game ends when there are no more moves available. The goal of Dominator is to minimize the number of vertices chosen, while Staller has the opposite goal. The Dominator-start domination game and the Staller-start domination game will be shortly referred to as the D-game and S-game, respectively. The {\em D-game domination number}, $\dstart(G)$, of $G$ is the number of moves in a D-game when both players play optimally. The {\em S-game domination number}, $\sstart(G)$, of $G$ is defined analogously for the S-game.

If $G$ is a graph and $S\subseteq V(G)$, then a partially dominated graph $G|S$ is a graph together with a declaration that the vertices from $S$ are already dominated. If $S = \{v\}$, then we will simplify the notation $G|\{v\}$ to $G|v$. We use $\dstart(G|S)$ (resp.\ $\sstart(G|S)$) to denote the number of moves remaining in the game on $G|S$ under optimal play when Dominator (resp.\ Staller) has the next move.

A graph $G$ is \emph{domination game critical}, {\em $\dstart$-critical} for short, if $\dstart(G|v) < \dstart(G)$ holds for every vertex $v \in V(G)$. If $G$ is $\dstart$-critical and $\dstart(G) = k$, we say that $G$ is $k$-$\dstart$-critical. This concept was introduced in~\cite{bujtas-2015c} where, among other results, a complete list of $\dstart$-critical trees on up to $17$ vertices was found by computer. The list consists of $13$ trees only. So it is natural to wonder if there exists an infinite family of $\dstart$-critical trees. Let $T_{p,q,r}$ be the tree obtained from disjoint paths $P_{4p+1}$, $P_{4q+1}$, and $P_{4r+1}$ by identifying three end-vertices, one from each path. It was conjectured in~\cite[p.~792]{bujtas-2015c} that the resulting three legged spider $T_{p,q,r}$ is $(2(p+q+r)+1)$-$\dstart$-critical. The conjecture has been verified by computer for $p,q,r\in \{0,1,2,3\}$ and it was stated that new techniques to prove the conjecture would be welcome. Using recent developments from~\cite{kosmrlj-2017} and additional techniques developed in this paper we can now confirm the conjecture.

The paper is organized as follows. In the next section we state concepts and recall some results needed in this paper. Then, in Section~\ref{sec:cutting-lemma}, the Cutting Lemma and the Union Lemma are proved. In the subsequent section we use these tools to prove the above conjecture that $T_{p,q,r}$ is $(2(p+q+r)+1)$-$\dstart$-critical.  In Section~\ref{sec:critical-trees} we report on our computational results and list all $\gamma_g$-critical trees on up to 20 vertices. It turned out that the variety of such trees is larger than earlier anticipated. We conclude the paper with a section in which the Cutting Lemma is extended such that it involves also the Staller-pass game.

\section{Preliminaries}
\label{sec:preliminaries}

A vertex of a graph \emph{dominates} itself and its neighbors; a \emph{dominating set} in a graph $G$ is a set of vertices of $G$ that dominates all vertices in the graph. We will say that a vertex selected according to the rules of the domination game is a {\em legal move}.
We next collect known results that will be useful to us.

\begin{thm}{\rm (\cite[Theorem~4.6]{kinnersley-2013})}
  \label{thm:tree-no-minus}
  If $F$ is a partially dominated forest, then $\dstart(F) \le \sstart(F)$.
\end{thm}

\begin{lemma}{\rm (\cite[Lemma~2.1]{kinnersley-2013})}
  \label{lem:continuation}
  {\rm (Continuation Principle)}
  Let $G$ be a graph, and let $A,B\subseteq V(G)$. If $B\subseteq A$, then
  $\dstart(G|A)\le \dstart(G|B)$ and $\sstart(G|A)\le \sstart(G|B)$.
\end{lemma}

We say that a graph $G$ {\em realizes} the pair $(k,\ell)$ if $\dstart(G)=k$ and $\sstart(G)=\ell$. The Continuation Principle in particular implies that $|\dstart(G) - \sstart(G)| \le 1$ holds for every graph $G$. Thus, every graph realizes either $(k,k+1)$, $(k,k)$, or $(k,k-1)$ for some integer~$k$. The graph is then consequently called a {\em plus graph}, an {\em equal graph}, or a {\em minus graph}, respectively. We also say that $G$ is a {\em no-minus} graph if $\dstart(G|S) \le \sstart(G|S)$ holds for any $S\subseteq V(G)$. (Note that Theorem~\ref{thm:tree-no-minus} asserts that forests are no-minus graphs.)

\begin{proposition}{\rm (\cite[Proposition~5]{bujtas-2015c})}
  \label{prp:no-plus}
  If $G$ is a $\dstart$-critical graph, then $G$ is either a minus graph or an equal graph.
\end{proposition}

\begin{thm}{\rm (\cite[Theorem~3]{bujtas-2015c})}
  \label{thm:critical}
  If $u$ is a vertex of a graph $G$, then $\dstart(G|u) \ge \dstart(G) - 2$ holds. Moreover, if $G$ is a no-minus graph, then $\dstart(G|u) \ge \dstart(G) - 1$.
\end{thm}

Following the notation from~\cite{kosmrlj-2017} we set $P_n' = P_{n+1}|u$ and $P_n'' = P_{n+2}|\{u,v\}$, where $u$ and $v$ are the vertices of the path of degree~$1$.

\begin{lemma}{\rm (\cite[Lemmas~2.1 and 2.3]{kosmrlj-2017})}
  \label{lem:P'-P''}
  If $n\ge 0$, then
  \begin{eqnarray*}
  \dstart(P'_n) = \dstart(P''_n)&\!\!\!=\!\!\!&\left\{\begin{array}{ll}\lceil \frac{n}{2}\rceil-1;& n \equiv 3 \, (\modo \, 4), \1 \\\lceil \frac{n}{2}\rceil,& {\rm otherwise,} \end{array}\right. \1\\
  \sstart(P'_n) = \sstart(P''_n)&\!\!\! = \!\!\! &\left\{\begin{array}{ll}\lceil \frac{n}{2}\rceil+1,& n\equiv 2 \, (\modo \, 4), \1 \\\lceil \frac{n}{2}\rceil;& {\rm otherwise.} \end{array}\right.
  \end{eqnarray*}
\end{lemma}

Underlying this result, there is a general observation that can be proved by a simple induction on the length of the path.
\begin{observation}\label{o:p'isp''}
For any graph $G$ and $n\ge0$, $\dstart(G\cup P'_n) = \dstart(G\cup P''_n)$ and $\sstart(G\cup P'_n) = \sstart(G \cup P''_n)$.
\end{observation}

A variation of the domination game when Dominator (resp.\ Staller) is allowed, but not obligated, to skip exactly one move in the course of the game, is called the {\em Dominator-pass game} (resp.\ {\em Staller-pass game}).  The number of moves in such a game, where both players are playing optimally, is denoted by $\gamma_{g}^{{\rm dp}}(G)$ (resp.\ $\gamma_{g}^{{\rm sp}}(G)$) when Dominator starts the game (unless he decides to pass already the first move) and by  $\gamma_{g}'^{{\rm dp}}(G)$ (resp.\ $\gamma_{g}'^{{\rm sp}}(G)$) when Staller starts the game.

\begin{proposition}{\rm (\cite[Lemma~2.2, Proposition~2.3]{dorbec-2015})}
  \label{lem:pass-in-general}
  If $S$ is a subset of vertices of a graph $G$, then $\gamma_g^{{\rm sp}}(G|S)\le\gamma_g(G|S)+1$. Moreover, if $G$ is a no-minus graph, then $\gamma_g^{{\rm sp}}(G|S)=\gamma_g^{{\rm dp}}(G|S)=\gamma_g(G|S)$.
\end{proposition}

\section{The Cutting Lemma and the Union Lemma}
\label{sec:cutting-lemma}

In this section we prove the new tools, the \emph{Cutting Lemma} and the \emph{Union Lemma}. In order to state the Cutting Lemma, we introduce some additional notation. Let $G$ be a graph and $uv$ an arbitrary edge in $G$. We define $G_{uv}$ as the graph obtained from $G - uv$ by adding two new vertices $u'$ and $v'$, adding the two edges $uv'$ and $vu'$, and declaring that both $u'$ and $v'$ are dominated.

\begin{thm}
{\rm (Cutting Lemma)}
\label{thm:cutting-lemma}
Let $G$ be a graph, and let $A,B \subseteq V(G)$ where $B\subseteq A$. If $uv$ is an edge of $G$, then
\[
\dstart(G|A) \le \dstart(G_{uv}|B) \hspace*{0.5cm} \mbox{and} \hspace*{0.5cm} \sstart(G|A) \le \sstart(G_{uv}|B).
\]
\end{thm}

\proof
We proceed by induction on $|V(G) \setminus A|$. If $|V(G) \setminus A| = 0$, then $\dstart(G|A) = 0$ and $\dstart(G_{uv}|B) \ge 0$, and $\sstart(G|A) = 0$ and $\sstart(G_{uv}|B) \ge 0$. Thus, the desired result holds. This establishes the base case. Suppose that $|V(G) \setminus A| \ge 1$, and so the game is not yet complete in $G|A$.

We first prove that $\dstart(G|A) \le \dstart(G_{uv}|B)$.
Let $x$ be an optimal first move for Dominator in $G_{uv}|B$.
We can assume that $x$ is not $u'$,
for $N_{G_{uv}}[u']\subseteq N_{G_{uv}}[v]$ and so $v$ is always a move at least as good as $u'$.
For similar reason, $x$ is not $v'$ either.
From $B\subseteq A$, we directly infer $B\cup N_{G_{uv}}[x] \subseteq A\cup N_G[x]\cup \{u',v'\}$.
On the one hand, if $N_G[x]\setminus A \neq \emptyset$
(that is, $x$ is a legal move in $G|A$), then we can apply induction and obtain
\[
  \begin{array}{rcl@{\qquad}l}
    \dstart(G_{uv}|B) & =   & 1+ \sstart\left(G_{uv}|(B\cup N_{G_{uv}}[x])\right) & \\
                           & \ge & 1+ \sstart(G|(A\cup N_G[x])) & \mbox{(by induction)}\\
                           & \ge & \dstart(G|A) & \\
\end{array}\]
as $x$ may not be an optimal move in $G|A$.

On the other hand, if $x$ is not a legal move in $G|A$,
Dominator may choose any other move $y$
and we still have $B\cup N_{G_{uv}}[x] \subseteq A\subseteq A\cup N_G[y]\cup \{u',v'\}$. So by induction,
\[
  \begin{array}{rcl@{\qquad}l}
    \dstart(G_{uv}|B) & =   & 1+ \sstart(G_{uv}|(B\cup N_{G_{uv}}[x])) & \\
                           & \ge & 1+ \sstart(G|(A\cup N_G[y])) & \mbox{(by induction)}\\
                           & \ge & \dstart(G|A) & \\
\end{array}\]
as $y$ may not be an optimal move in $G|A$.

Let us now prove that $\sstart(G|A) \le \sstart(G_{uv}|B)$.
Consider an optimal move $x$ for Staller in $G|A$.
We have that $N_{G_{uv}}[x]\setminus\{u',v'\}\subseteq N_G[x]$, hence since $B\subseteq A$, we also have $B\cup N_{G_{uv}}[x] \subseteq A\cup N_G[x]\cup \{u',v'\}$.
On the one hand, if $N_{G_{uv}}[x]\setminus B$ is non empty, then
\[
  \begin{array}{rcl@{\qquad}l}
  \sstart(G|A) & =   & 1+ \dstart(G|(A\cup N_G[x])) & \mbox{(}x\mbox{ is optimal for Staller)} \\
                  & \le & 1+ \dstart(G_{uv}|(B\cup N_{G_{uv}}[x])) & \mbox{(by induction)}\\
                  & \le & \sstart(G_{uv}|B) & \\
\end{array}\]
since $x$ is not necessarily an optimal move for Staller in $G_{uv}|B$.

On the other hand, if $N_{G_{uv}}[x]\setminus B$ is empty, that means that $x$ is $u$ or $v$ and the only newly dominated vertex in $G|A$ is $v$ or $u$.
Renaming vertices if necessary, suppose $x=u$.
Then $u'$ is a legal move in $G_{uv}|B$ that newly dominates exactly the same vertex, that is $B\cup N_{G_{uv}}[u']\subseteq A \cup N_G[x]\cup \{u',v'\}$ and we apply induction with
\[
  \begin{array}{rcl@{\qquad}l}
  \sstart(G|A) & =   & 1+ \dstart(G|(A\cup N_G[x])) & \mbox{(}x\mbox{ is optimal for Staller)} \\
                  & \le & 1+ \dstart(G_{uv}|(B\cup N_{G_{uv}}[u'])) & \mbox{(by induction)}\\
                  & \le & \sstart(G_{uv}|B) & \\
\end{array}\]
since $u'$ may not be an optimal move for Staller in $G_{uv}|B$.
That concludes the proof.
\qed

It may be worth noting that the above proof differs from the proofs by imagination strategy only by its presentation. In the first half, $G_{uv}|B$ plays the role of the imagined game for Dominator, and in the second, $G|A$ does so for Staller.

Recall that $P_n' = P_{n+1}|u$ and $P_n'' = P_{n+2}|\{u,v\}$, where $u$ and $v$ are the vertices of the path of degree~$1$. Further recall that by Lemma~\ref{lem:P'-P''}, $\sstart(P'_n) = \sstart(P''_n)$.  Defining the {\em weighting function} $\w$ of paths with
\[
\w(P'_{4q+r})= \w(P''_{4q+r})= 2q +
  \begin{cases}
    0 \mbox{ if } r=0\\
    1  \mbox{ if } r=1\\
    \frac{3}{2} \mbox{ if } r=2\1 \\
    \frac{7}{4} \mbox{ if } r=3,
  \end{cases}
\]
the Union Lemma reads as follows.

\begin{lemma}[Union Lemma]
\label{lem:union}
If $F_1, \ldots, F_k$ are vertex disjoint paths where $F_i = P_{n_i}'$ or $F_i = P_{n_i}''$ for $i \in [k]$ and $n_i \ge 1$, then
\[
\sstart\left(\bigcup_{i = 1}^k F_i \right) \le \left\lceil \sum_{i = 1}^k \w(F_i) \right\rceil.
\]
\end{lemma}

\proof
Let $F$ be the union of the vertex disjoint paths $F_1, \ldots, F_k$ where $F_i = P_{n_i}'$ or $F_i = P_{n_i}''$ for $i \in [k]$ and $n_i \ge 1$, and so
\[
F = \bigcup_{i = 1}^k F_i.
\]
Further, let denote
\[
\varphi(F) = \sum_{i = 1}^k \w(F_i).
\]

We proceed by induction on $k \ge 1$ to show that $\sstart(F) \le \lceil \varphi(F)  \rceil$. If $k = 1$, then the desired result follows from Lemma~\ref{lem:P'-P''}. This establishes the base case. Suppose that $k \ge 2$ and that the desired results holds for the disjoint union of fewer than $k$ paths.

We first observe that by Observation~\ref{o:p'isp''} and the Cutting Lemma, we have for a general graph $G$ that $\sstart(G \cup P'_{n+4}) = \sstart(G \cup P''_{n+4}) \le \sstart(G \cup P''_4 \cup P_n'')$. Splitting long paths in such a way, we may therefore assume that every path in the union $F$ has at least one but at most four undominated (internal) vertices with both ends of the path declared dominated.

Let $F'$ be the disjoint union of paths that results after Staller and Dominator have played their first moves.  To prove the union lemma, we simply describe Dominator's strategy to respond to each of Staller moves, so that what remains after Dominator's move is a union $F'$ of paths such that $\left\lceil\varphi(F)\right\rceil\ge \left\lceil \varphi(F')\right\rceil+2$. With this, the desired result follows readily by induction.

\smallskip
\noindent \textbf{Case 1.} Suppose Staller dominates some new vertex of a path $P_4''$ in the union $F$.   Whichever move Staller made, Dominator dominates in his move all the yet undominated vertices of this path. Thus, $F = F' \cup P''_4$ and $\varphi(F)=\varphi(F')+\w(P''_4) = \varphi(F')+2$.

\smallskip
\noindent \textbf{Case 2.}
Suppose Staller's move dominates some vertex of a path $P''_3$. This path then becomes a $P''_2$ or smaller. If there exists another $P''_3$ in $F$, then Dominator plays as his first move its middle vertex, dominating it in one turn, and we get $\varphi(F) \ge   \varphi(F') - \w(P''_2)+2\w(P''_3) = \varphi(F')+2$.  Otherwise, Dominator finishes dominating the same $P''_3$ as Staller. In that case, it means that $\varphi(F')= \frac{\ell}{2}$ for some integer $\ell$, since in the case for all $P \in F'$, $\w(P)$ is a multiple of $\frac{1}{2}$. Thus,
\[
\left\lceil{\varphi(F)}\right\rceil =\left\lceil \varphi(F')+\frac{7}{4}\right\rceil=\left\lceil \frac{\ell}{2}+\frac{7}{4}\right\rceil
  =\left\lceil \frac{\ell}{2}+\frac{7}{4}+\frac{1}{4}\right\rceil = \left\lceil\frac{\ell}{2}\right\rceil+2=
  \left\lceil\varphi(F')\right\rceil+2.
\]

\smallskip
\noindent \textbf{Case 3.}
Suppose now Staller's move dominates some vertex of a path $P''_2$, leaving at most one undominated vertex. If there exists another $P''_2$ or $P''_3$ in $F$,  Dominator dominates one of them with his move, and $\varphi(F)\ge \varphi(F')-\w(P''_1)+\w(P''_2)+ \w(P) \ge \varphi(F')+2$, where $P \in \{P_2', P_3'\}$. Otherwise, Dominator dominates the end of the same path $P''_2$ Staller played on. In that case, $\varphi(F')$ is an integer, and thus
\[
\left\lceil{\varphi(F)}\right\rceil
= \left\lceil \varphi(F')+\frac{3}{2}\right\rceil
= \varphi(F')+\left\lceil\frac{3}{2}\right\rceil
= \left\lceil\varphi(F')\right\rceil+2.
\]

\smallskip
\noindent \textbf{Case 4.}
If Staller's move dominates the only undominated vertex of a path $P''_1$, then Dominator either dominates any path with three or less undominated vertices in one move, or dominates three vertices from a $P''_4$, leaving it a $P''_1$. In either of the cases we have $\varphi(F)\ge \varphi(F')+2$ and we are done.
\qed

\section{Three legged spiders are $\dstart$-critical}
\label{sec:spiders}

Using the tools developed in the previous section we now prove the following result conjectured in~\cite{bujtas-2015c}.

\begin{thm}
\label{thm:main}
If $p,q,r\ge 1$, then $T_{p,q,r}$ is a $(2(p+q+r)+1)$-$\dstart$-critical graph.
\end{thm}

\proof
We first prove that $\dstart(T_{p,q,r})\ge 2(p+q+r)+1$.
  For this we describe a strategy for Staller to ensure at least $2(p+q+r)+1$ moves are played.
  After each move of Dominator, she tries to dominate only one new vertex,
  and when possible in the same leg that Dominator just played.
  Since Dominator has the first move,
  she can always play a vertex already dominated
  and dominate one new vertex if that vertex is of degree 2 or maybe two if the only legal move of Staller on a dominated vertex is the vertex of degree~$3$. The latter case may happen only once.
  Note also that Dominator may dominate four vertices in one move only once by playing the vertex of degree~$3$,
  otherwise he may not dominate more than three new vertices at a time.

\medskip \noindent \textbf{Case 1.}
Suppose first that Dominator dominates at most three new vertices on each of his moves
  and Staller only one new vertex on each of her moves.
  Then, after Dominator played his last move, at most one vertex remains undominated.
  Denoting with $m$ the number of moves of Dominator, we have $3m+(m-1)\ge 4p+4q+4r$ and thus $m\ge p+q+r+1$ (since it is an integer).
  Thus the total number of moves in the game is at least $m+m-1 = 2(p+q+r)+1$,
  and Staller's strategy succeeded. Note that the same computation works also if Dominator plays once a move that dominates four vertices but
  also plays a move that dominates two or less vertices.

\medskip
\noindent \textbf{Case 2.}
  Suppose next that Dominator forces Staller to dominate two vertices in one move.
  That is, Staller is forced to play the vertex of degree~$3$ and doing that to cover two new vertices.
  This means that the previously dominated vertices are exactly all the vertices of a branch. Assume without loss of generality that the $4p+1$ vertices of the $P_{4p+1}$ leg were dominated after Dominator's move.
  Letting $m_1$ be the total number of moves played by Dominator before that move,
  Staller necessarily played $m'_1=m_1-1$ moves.
  By our assumption that Dominator dominates at most three new vertices on each of his moves and that Staller
  dominates one new vertex on each of her $m_1-1$ first moves, we get that
  $3 m_1 + (m_1-1) \ge 4p+1$. Since $m_1$ is an integer, we have $m_1\ge p+1$.
  Then Staller dominates two new vertices, and there are $4q-1+4r-1$ vertices left to dominate.
  Let $m_2$ be the number of moves Dominator plays in the rest of the game.
  After Dominator plays his last move, at most one vertex is not yet dominated,
  and Staller has played $m_2-1$ additional moves.
  So $3m_2+(m_2-1)\ge 4q+4r-3$  and thus $m_2\ge q+r$.
  Thus, the total number of moves played in the course of the game is at least
  $2m_1+2m_2-1\ge 2p+2q+2r+1$, so Staller's strategy succeeded also in this case.

\medskip
\noindent \textbf{Case 3.}
  Suppose now that Dominator used the opportunity to dominate four new vertices at some stage of the game by playing the vertex of degree~$3$.
  Note that in the earlier course of the game, after each move of Dominator,
  Staller was able to dominate one new vertex in the same branch,
  and that the number of dominated vertices in each leg is thus a multiple of~$4$.
  So the number of undominated vertices in each leg after Dominator played the center vertex is $3\bmod{4}$.
  We now need to slightly refine Staller's strategy.

  We call the \emph{state of the game} the number of undominated vertices $\bmod~4$ in each branch, described with a corresponding multiset. So the state of the game just after the move of Dominator on the center is $\{3,3,3\}$.
  A move of Staller brings the game into the state $\{3,3,2\}$.
  Now Dominator, while dominating three new vertices, can bring the game back to the state $\{3,3,3\}$
  while the number of undominated vertices in the corresponding leg is large enough, and Staller responds similarly.
  Eventually, Dominator is forced to play in another leg if he wants to dominate three new vertices,
  bringing the game to the state $\{3,0,2\}$.
  Then Staller should use the opportunity to bring the game to the state $\{2,0,2\}$.
  There Dominator can for some time push the game back to the state $\{3,0,2\}$ or to the state $\{2,1,2\}$,
  from where Staller again replies is such a way that the game returns to the state $\{2,0,2\}$.
  But eventually, Dominator is forced to finish a leg with only two undominated vertices left. Now  the same computation as done in Case~1 concludes the argument in this case.

From the above three cases we conclude that
\begin{equation}
\label{eq:whatever1}
\dstart(T_{p,q,r})\ge 2(p+q+r) + 1\,.
\end{equation}

We next prove that $\dstart(T_{p,q,r}|v)\le 2(p+q+r)$ holds for every vertex $v$. Without loss of generality, assume that $v$ belongs to the leg of length $4p$. Denote the vertices of this leg $v_0,v_1,\ldots,v_{4p}$ where $v_{4p}=c$ is the vertex of degree~$3$ in the spider.

Assume first that $v$ is some $v_i$ with $i\not\equiv 1 \bmod{4}$ and $i\le 4p-2$. Let Dominator's first move be on $c$.
  We use the Cutting Lemma on the edge $v_iv_{i+1}$ which isolates a $P'_{i}$,
  and on the edges linking $c$ to the paths, obtaining a path $P'_{4p-i-2}$, a path $P'_{4q-1}$ and a path $P'_{4r-1}$  such that $\dstart(T_{p,q,r}|v) \le 1+\sstart\big(\bigcup(P'_i,P'_{4p-i-2},P'_{4q-1},P'_{4r-1})\big)$.
  We now apply the Union Lemma and get:
  \begin{itemize}
    \item
      If $i\equiv 0$ or $2 \bmod{4}$, then
      \[
        \dstart(T_{p,q,r}|v) \le 1 +
        \left\lceil2(p-1) + \frac{3}{2} + 2(q-1) + \frac{7}{4} + 2(r-1) + \frac{7}{4}\right\rceil = 2(p+q+r)\,.
      \]
    \item
      If $i\equiv 3 \bmod{4}$, then
      \[
        \dstart(T_{p,q,r}|v) \le 1 +
        \left\lceil2(p-2) + \frac{7}{4} + \frac{7}{4} + 2(q-1) + \frac{7}{4} + 2(r-1) + \frac{7}{4}\right\rceil = 2(p+q+r)\,.
      \]
  \end{itemize}

  Suppose now that $i\equiv 1 \bmod{4}$. Dominator's response is then on $v_{4p-1}$.
  We use the Cutting Lemma on the edges $v_iv_{i+1}$ which isolates a $P'_i$, and on the edges linking $c$ to the paths,
  obtaining path $P'_{4p-i-3}$, a $P'_{4q}$ and a $P'_{4r}$.
Applying the Union Lemma again, we get
\[
\begin{array}{lcl}
\dstart(T_{p,q,r}|v) & \le & \displaystyle{ 1+\sstart\left(\bigcup(P'_i,P'_{4p-i-3},P'_{4q},P'_{4r})\right) } \1 \\
& \le & 1+ 2(p-1)+ 1 +2q +2r \1 \\
& = & 2(p+q+r).
\end{array}
\]

There are two cases left, when $i$ is $4p$ and $4p-1$. In both cases, the first move of Dominator is $v_{i-2}$.
  We then use the Cutting Lemma on the edge between $c$ and the rest of the length $4q$-path (whose both ends are undominated in the second case),
  and on the edge between $v_{i-4}$ and $v_{i-3}$.
  If $i=4p$, then we get the union of the paths $P'_{4p-3}, P'_{4q}, P'_{4r}$ and thus
  \[ \dstart(T_{p,q,r}|v) \le 1 +
    \left\lceil2(p-1) + 1 + 2q + 2r \right\rceil = 2(p+q+r)\,.
  \]
  Otherwise, if $i=4p-1$, then we get the union of the paths $P'_{4p-4}$, $P'_{4q}$ and $P'_{4r+1}$. We thus have
  \[ \dstart(T_{p,q,r}|v) \le 1 +
    \left\lceil2(p-1) + 2q + 2r + 1 \right\rceil = 2(p+q+r)\,.
  \]
  Thus for every vertex $v$,
\begin{equation}
\label{eq:whatever2}
\dstart(T_{p,q,r}|v)\le 2(p+q+r).
\end{equation}
  Now, let $v$ be an optimal first move for Staller in the S-game on $T_{p,q,r}$. We have
$$
\begin{array}{lcll}
2(p+q+r)+1 & \le & \dstart(T_{p,q,r}) & (\textrm{by } \eqref{eq:whatever1}) \\
& \le & \sstart(T_{p,q,r}) & (\textrm{by Theorem}~\ref{thm:tree-no-minus})\\
& \le & 1+ \dstart(T_{p,q,r}|N[v]) & (v \textrm{ is optimal for Staller}) \\
& \le & 1+ \dstart(T_{p,q,r}|v) & (\textrm{Continuation Principle})\\
& \le & 1+2(p+q+r)  & (\textrm{by } \eqref{eq:whatever2})\,.
\end{array}
$$
So we have equality throughout and in particular, $\dstart(T_{p,q,r})= 2(p+q+r)+1$.
\qed

Note that the last chain of inequalities from the above proof implies that
$$\sstart(T_{p,q,r}) = 2(p+q+r) + 1$$
holds for $p,q,r\ge 1$. This fact could also be deduced as follows. By Theorem~\ref{thm:tree-no-minus} we know that $T_{p,q,r}$ is either equal or plus. Since by Theorem~\ref{thm:main}, $T_{p,q,r}$ is $\dstart$-critical, Proposition~\ref{prp:no-plus} implies that $T_{p,q,r}$ cannot be plus. Hence it is equal.

\section{The variety of $\gamma_g$-critical trees}
\label{sec:critical-trees}

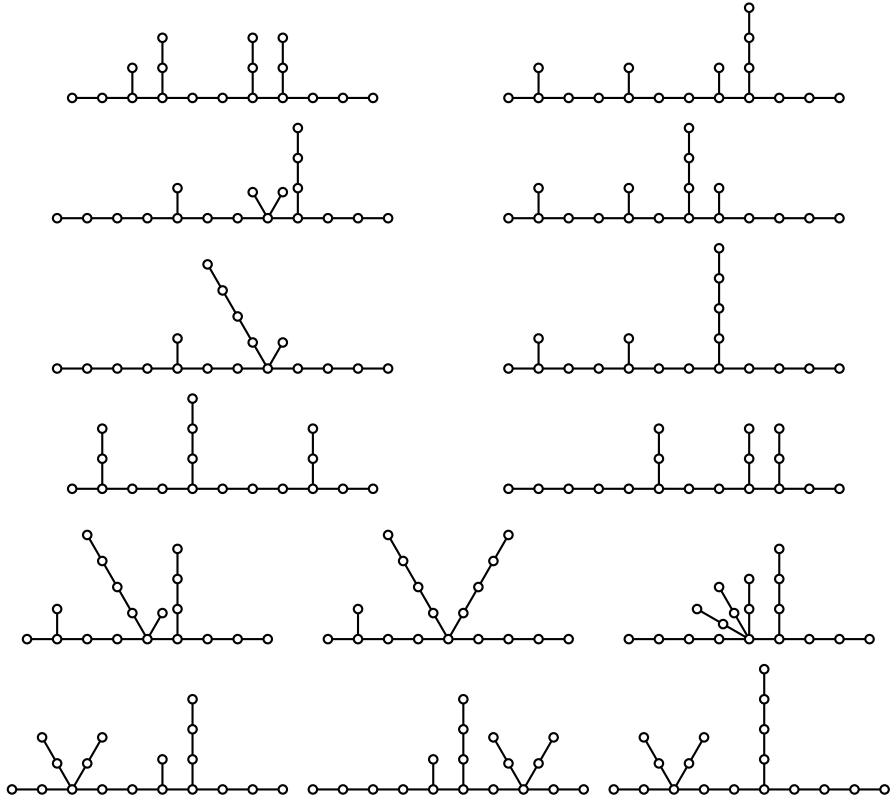
\begin{figure}[ht!]
\begin{center}
\begin{tikzpicture}[scale=0.4,style=thick]
\def\vr{2pt/0.5}

\coordinate(g18x1y1) at (2.0,-4);
\begin{scope}[shift=(g18x1y1)]
\foreach \i in {0,...,10}
{
	\coordinate(x\i) at (\i,0);
}

\draw(x0) -- (x1) -- (x2) -- (x3) -- (x4) -- (x5) -- (x6) -- (x7) -- (x8) -- (x9) -- (x10);

\coordinate(x11) at ({2+1*cos(90.0)}, {1*sin(90.0)});
\draw(x2) -- (x11);

\coordinate(x12) at ({3+1*cos(90.0)}, {1*sin(90.0)});
\coordinate(x13) at ({3+2*cos(90.0)}, {2*sin(90.0)});
\draw(x3) -- (x12) -- (x13);

\coordinate(x14) at ({6+1*cos(90.0)}, {1*sin(90.0)});
\coordinate(x15) at ({6+2*cos(90.0)}, {2*sin(90.0)});
\draw(x6) -- (x14) -- (x15);

\coordinate(x16) at ({7+1*cos(90.0)}, {1*sin(90.0)});
\coordinate(x17) at ({7+2*cos(90.0)}, {2*sin(90.0)});
\draw(x7) -- (x16) -- (x17);

\foreach \i in {0,...,17}
{
	\draw (x\i)[fill=white] circle (\vr);
}
\end{scope}

\coordinate(g18x1y2) at (16.5,-4);
\begin{scope}[shift=(g18x1y2)]
\foreach \i in {0,...,11}
{
	\coordinate(x\i) at (\i,0);
}

\draw(x0) -- (x1) -- (x2) -- (x3) -- (x4) -- (x5) -- (x6) -- (x7) -- (x8) -- (x9) -- (x10) -- (x11);

\coordinate(x12) at ({1+1*cos(90.0)}, {1*sin(90.0)});
\draw(x1) -- (x12);

\coordinate(x13) at ({4+1*cos(90.0)}, {1*sin(90.0)});
\draw(x4) -- (x13);

\coordinate(x14) at ({7+1*cos(90.0)}, {1*sin(90.0)});
\draw(x7) -- (x14);

\coordinate(x15) at ({8+1*cos(90.0)}, {1*sin(90.0)});
\coordinate(x16) at ({8+2*cos(90.0)}, {2*sin(90.0)});
\coordinate(x17) at ({8+3*cos(90.0)}, {3*sin(90.0)});
\draw(x8) -- (x15) -- (x16) -- (x17);

\foreach \i in {0,...,17}
{
	\draw (x\i)[fill=white] circle (\vr);
}
\end{scope}

\coordinate(g18x2y1) at (1.5,-8);
\begin{scope}[shift=(g18x2y1)]
\foreach \i in {0,...,11}
{
	\coordinate(x\i) at (\i,0);
}

\draw(x0) -- (x1) -- (x2) -- (x3) -- (x4) -- (x5) -- (x6) -- (x7) -- (x8) -- (x9) -- (x10) -- (x11);

\coordinate(x12) at ({4+1*cos(90.0)}, {1*sin(90.0)});
\draw(x4) -- (x12);

\coordinate(x13) at ({7+1*cos(60.0)}, {1*sin(60.0)});
\draw(x7) -- (x13);

\coordinate(x14) at ({7+1*cos(120.0)}, {1*sin(120.0)});
\draw(x7) -- (x14);

\coordinate(x15) at ({8+1*cos(90.0)}, {1*sin(90.0)});
\coordinate(x16) at ({8+2*cos(90.0)}, {2*sin(90.0)});
\coordinate(x17) at ({8+3*cos(90.0)}, {3*sin(90.0)});
\draw(x8) -- (x15) -- (x16) -- (x17);

\foreach \i in {0,...,17}
{
	\draw (x\i)[fill=white] circle (\vr);
}
\end{scope}

\coordinate(g18x2y2) at (16.5,-8);
\begin{scope}[shift=(g18x2y2)]
\foreach \i in {0,...,11}
{
	\coordinate(x\i) at (\i,0);
}

\draw(x0) -- (x1) -- (x2) -- (x3) -- (x4) -- (x5) -- (x6) -- (x7) -- (x8) -- (x9) -- (x10) -- (x11);

\coordinate(x12) at ({1+1*cos(90.0)}, {1*sin(90.0)});
\draw(x1) -- (x12);

\coordinate(x13) at ({4+1*cos(90.0)}, {1*sin(90.0)});
\draw(x4) -- (x13);

\coordinate(x14) at ({6+1*cos(90.0)}, {1*sin(90.0)});
\coordinate(x15) at ({6+2*cos(90.0)}, {2*sin(90.0)});
\coordinate(x16) at ({6+3*cos(90.0)}, {3*sin(90.0)});
\draw(x6) -- (x14) -- (x15) -- (x16);

\coordinate(x17) at ({7+1*cos(90.0)}, {1*sin(90.0)});
\draw(x7) -- (x17);

\foreach \i in {0,...,17}
{
	\draw (x\i)[fill=white] circle (\vr);
}
\end{scope}

\coordinate(g18x3y1) at (1.5,-13);
\begin{scope}[shift=(g18x3y1)]
\foreach \i in {0,...,11}
{
	\coordinate(x\i) at (\i,0);
}

\draw(x0) -- (x1) -- (x2) -- (x3) -- (x4) -- (x5) -- (x6) -- (x7) -- (x8) -- (x9) -- (x10) -- (x11);

\coordinate(x12) at ({4+1*cos(90.0)}, {1*sin(90.0)});
\draw(x4) -- (x12);

\coordinate(x13) at ({7+1*cos(60.0)}, {1*sin(60.0)});
\draw(x7) -- (x13);

\coordinate(x14) at ({7+1*cos(120.0)}, {1*sin(120.0)});
\coordinate(x15) at ({7+2*cos(120.0)}, {2*sin(120.0)});
\coordinate(x16) at ({7+3*cos(120.0)}, {3*sin(120.0)});
\coordinate(x17) at ({7+4*cos(120.0)}, {4*sin(120.0)});
\draw(x7) -- (x14) -- (x15) -- (x16) -- (x17);

\foreach \i in {0,...,17}
{
	\draw (x\i)[fill=white] circle (\vr);
}
\end{scope}

\coordinate(g18x3y2) at (16.5,-13);
\begin{scope}[shift=(g18x3y2)]
\foreach \i in {0,...,11}
{
	\coordinate(x\i) at (\i,0);
}

\draw(x0) -- (x1) -- (x2) -- (x3) -- (x4) -- (x5) -- (x6) -- (x7) -- (x8) -- (x9) -- (x10) -- (x11);

\coordinate(x12) at ({1+1*cos(90.0)}, {1*sin(90.0)});
\draw(x1) -- (x12);

\coordinate(x13) at ({4+1*cos(90.0)}, {1*sin(90.0)});
\draw(x4) -- (x13);

\coordinate(x14) at ({7+1*cos(90.0)}, {1*sin(90.0)});
\coordinate(x15) at ({7+2*cos(90.0)}, {2*sin(90.0)});
\coordinate(x16) at ({7+3*cos(90.0)}, {3*sin(90.0)});
\coordinate(x17) at ({7+4*cos(90.0)}, {4*sin(90.0)});
\draw(x7) -- (x14) -- (x15) -- (x16) -- (x17);

\foreach \i in {0,...,17}
{
	\draw (x\i)[fill=white] circle (\vr);
}
\end{scope}

\coordinate(g18x4y1) at (2.0,-17);
\begin{scope}[shift=(g18x4y1)]
\foreach \i in {0,...,10}
{
	\coordinate(x\i) at (\i,0);
}

\draw(x0) -- (x1) -- (x2) -- (x3) -- (x4) -- (x5) -- (x6) -- (x7) -- (x8) -- (x9) -- (x10);

\coordinate(x11) at ({1+1*cos(90.0)}, {1*sin(90.0)});
\coordinate(x12) at ({1+2*cos(90.0)}, {2*sin(90.0)});
\draw(x1) -- (x11) -- (x12);

\coordinate(x13) at ({4+1*cos(90.0)}, {1*sin(90.0)});
\coordinate(x14) at ({4+2*cos(90.0)}, {2*sin(90.0)});
\coordinate(x15) at ({4+3*cos(90.0)}, {3*sin(90.0)});
\draw(x4) -- (x13) -- (x14) -- (x15);

\coordinate(x16) at ({8+1*cos(90.0)}, {1*sin(90.0)});
\coordinate(x17) at ({8+2*cos(90.0)}, {2*sin(90.0)});
\draw(x8) -- (x16) -- (x17);

\foreach \i in {0,...,17}
{
	\draw (x\i)[fill=white] circle (\vr);
}
\end{scope}

\coordinate(g18x4y2) at (16.5,-17);
\begin{scope}[shift=(g18x4y2)]
\foreach \i in {0,...,11}
{
	\coordinate(x\i) at (\i,0);
}

\draw(x0) -- (x1) -- (x2) -- (x3) -- (x4) -- (x5) -- (x6) -- (x7) -- (x8) -- (x9) -- (x10) -- (x11);

\coordinate(x12) at ({5+1*cos(90.0)}, {1*sin(90.0)});
\coordinate(x13) at ({5+2*cos(90.0)}, {2*sin(90.0)});
\draw(x5) -- (x12) -- (x13);

\coordinate(x14) at ({8+1*cos(90.0)}, {1*sin(90.0)});
\coordinate(x15) at ({8+2*cos(90.0)}, {2*sin(90.0)});
\draw(x8) -- (x14) -- (x15);

\coordinate(x16) at ({9+1*cos(90.0)}, {1*sin(90.0)});
\coordinate(x17) at ({9+2*cos(90.0)}, {2*sin(90.0)});
\draw(x9) -- (x16) -- (x17);

\foreach \i in {0,...,17}
{
	\draw (x\i)[fill=white] circle (\vr);
}
\end{scope}

\coordinate(g18x5y1) at (0.5,-22);
\begin{scope}[shift=(g18x5y1)]
\foreach \i in {0,...,8}
{
	\coordinate(x\i) at (\i,0);
}

\draw(x0) -- (x1) -- (x2) -- (x3) -- (x4) -- (x5) -- (x6) -- (x7) -- (x8);

\coordinate(x9) at ({1+1*cos(90.0)}, {1*sin(90.0)});
\draw(x1) -- (x9);

\coordinate(x10) at ({4+1*cos(60.0)}, {1*sin(60.0)});
\draw(x4) -- (x10);

\coordinate(x11) at ({4+1*cos(120.0)}, {1*sin(120.0)});
\coordinate(x12) at ({4+2*cos(120.0)}, {2*sin(120.0)});
\coordinate(x13) at ({4+3*cos(120.0)}, {3*sin(120.0)});
\coordinate(x14) at ({4+4*cos(120.0)}, {4*sin(120.0)});
\draw(x4) -- (x11) -- (x12) -- (x13) -- (x14);

\coordinate(x15) at ({5+1*cos(90.0)}, {1*sin(90.0)});
\coordinate(x16) at ({5+2*cos(90.0)}, {2*sin(90.0)});
\coordinate(x17) at ({5+3*cos(90.0)}, {3*sin(90.0)});
\draw(x5) -- (x15) -- (x16) -- (x17);

\foreach \i in {0,...,17}
{
	\draw (x\i)[fill=white] circle (\vr);
}
\end{scope}

\coordinate(g18x5y2) at (10.5,-22);
\begin{scope}[shift=(g18x5y2)]
\foreach \i in {0,...,8}
{
	\coordinate(x\i) at (\i,0);
}

\draw(x0) -- (x1) -- (x2) -- (x3) -- (x4) -- (x5) -- (x6) -- (x7) -- (x8);

\coordinate(x9) at ({1+1*cos(90.0)}, {1*sin(90.0)});
\draw(x1) -- (x9);

\coordinate(x10) at ({4+1*cos(60.0)}, {1*sin(60.0)});
\coordinate(x11) at ({4+2*cos(60.0)}, {2*sin(60.0)});
\coordinate(x12) at ({4+3*cos(60.0)}, {3*sin(60.0)});
\coordinate(x13) at ({4+4*cos(60.0)}, {4*sin(60.0)});
\draw(x4) -- (x10) -- (x11) -- (x12) -- (x13);

\coordinate(x14) at ({4+1*cos(120.0)}, {1*sin(120.0)});
\coordinate(x15) at ({4+2*cos(120.0)}, {2*sin(120.0)});
\coordinate(x16) at ({4+3*cos(120.0)}, {3*sin(120.0)});
\coordinate(x17) at ({4+4*cos(120.0)}, {4*sin(120.0)});
\draw(x4) -- (x14) -- (x15) -- (x16) -- (x17);

\foreach \i in {0,...,17}
{
	\draw (x\i)[fill=white] circle (\vr);
}
\end{scope}

\coordinate(g18x5y3) at (20.5,-22);
\begin{scope}[shift=(g18x5y3)]
\foreach \i in {0,...,8}
{
	\coordinate(x\i) at (\i,0);
}

\draw(x0) -- (x1) -- (x2) -- (x3) -- (x4) -- (x5) -- (x6) -- (x7) -- (x8);

\coordinate(x9) at ({4+1*cos(120.0)}, {1*sin(120.0)});
\coordinate(x10) at ({4+2*cos(120.0)}, {2*sin(120.0)});
\draw(x4) -- (x9) -- (x10);

\coordinate(x11) at ({4+1*cos(90.0)}, {1*sin(90.0)});
\coordinate(x12) at ({4+2*cos(90.0)}, {2*sin(90.0)});
\draw(x4) -- (x11) -- (x12);

\coordinate(x13) at ({4+1*cos(150.0)}, {1*sin(150.0)});
\coordinate(x14) at ({4+2*cos(150.0)}, {2*sin(150.0)});
\draw(x4) -- (x13) -- (x14);

\coordinate(x15) at ({5+1*cos(90.0)}, {1*sin(90.0)});
\coordinate(x16) at ({5+2*cos(90.0)}, {2*sin(90.0)});
\coordinate(x17) at ({5+3*cos(90.0)}, {3*sin(90.0)});
\draw(x5) -- (x15) -- (x16) -- (x17);

\foreach \i in {0,...,17}
{
	\draw (x\i)[fill=white] circle (\vr);
}
\end{scope}

\coordinate(g18x6y1) at (0.0,-27);
\begin{scope}[shift=(g18x6y1)]
\foreach \i in {0,...,9}
{
	\coordinate(x\i) at (\i,0);
}

\draw(x0) -- (x1) -- (x2) -- (x3) -- (x4) -- (x5) -- (x6) -- (x7) -- (x8) -- (x9);

\coordinate(x10) at ({2+1*cos(60.0)}, {1*sin(60.0)});
\coordinate(x11) at ({2+2*cos(60.0)}, {2*sin(60.0)});
\draw(x2) -- (x10) -- (x11);

\coordinate(x12) at ({2+1*cos(120.0)}, {1*sin(120.0)});
\coordinate(x13) at ({2+2*cos(120.0)}, {2*sin(120.0)});
\draw(x2) -- (x12) -- (x13);

\coordinate(x14) at ({5+1*cos(90.0)}, {1*sin(90.0)});
\draw(x5) -- (x14);

\coordinate(x15) at ({6+1*cos(90.0)}, {1*sin(90.0)});
\coordinate(x16) at ({6+2*cos(90.0)}, {2*sin(90.0)});
\coordinate(x17) at ({6+3*cos(90.0)}, {3*sin(90.0)});
\draw(x6) -- (x15) -- (x16) -- (x17);

\foreach \i in {0,...,17}
{
	\draw (x\i)[fill=white] circle (\vr);
}
\end{scope}

\coordinate(g18x6y2) at (10.0,-27);
\begin{scope}[shift=(g18x6y2)]
\foreach \i in {0,...,9}
{
	\coordinate(x\i) at (\i,0);
}

\draw(x0) -- (x1) -- (x2) -- (x3) -- (x4) -- (x5) -- (x6) -- (x7) -- (x8) -- (x9);

\coordinate(x10) at ({4+1*cos(90.0)}, {1*sin(90.0)});
\draw(x4) -- (x10);

\coordinate(x11) at ({5+1*cos(90.0)}, {1*sin(90.0)});
\coordinate(x12) at ({5+2*cos(90.0)}, {2*sin(90.0)});
\coordinate(x13) at ({5+3*cos(90.0)}, {3*sin(90.0)});
\draw(x5) -- (x11) -- (x12) -- (x13);

\coordinate(x14) at ({7+1*cos(60.0)}, {1*sin(60.0)});
\coordinate(x15) at ({7+2*cos(60.0)}, {2*sin(60.0)});
\draw(x7) -- (x14) -- (x15);

\coordinate(x16) at ({7+1*cos(120.0)}, {1*sin(120.0)});
\coordinate(x17) at ({7+2*cos(120.0)}, {2*sin(120.0)});
\draw(x7) -- (x16) -- (x17);

\foreach \i in {0,...,17}
{
	\draw (x\i)[fill=white] circle (\vr);
}
\end{scope}

\coordinate(g18x6y3) at (20.0,-27);
\begin{scope}[shift=(g18x6y3)]
\foreach \i in {0,...,9}
{
	\coordinate(x\i) at (\i,0);
}

\draw(x0) -- (x1) -- (x2) -- (x3) -- (x4) -- (x5) -- (x6) -- (x7) -- (x8) -- (x9);

\coordinate(x10) at ({2+1*cos(60.0)}, {1*sin(60.0)});
\coordinate(x11) at ({2+2*cos(60.0)}, {2*sin(60.0)});
\draw(x2) -- (x10) -- (x11);

\coordinate(x12) at ({2+1*cos(120.0)}, {1*sin(120.0)});
\coordinate(x13) at ({2+2*cos(120.0)}, {2*sin(120.0)});
\draw(x2) -- (x12) -- (x13);

\coordinate(x14) at ({5+1*cos(90.0)}, {1*sin(90.0)});
\coordinate(x15) at ({5+2*cos(90.0)}, {2*sin(90.0)});
\coordinate(x16) at ({5+3*cos(90.0)}, {3*sin(90.0)});
\coordinate(x17) at ({5+4*cos(90.0)}, {4*sin(90.0)});
\draw(x5) -- (x14) -- (x15) -- (x16) -- (x17);

\foreach \i in {0,...,17}
{
	\draw (x\i)[fill=white] circle (\vr);
}
\end{scope}

\end{tikzpicture}
\end{center}
\caption{All critical trees on 18 vertices}
\label{fig:critical-trees-18}
\end{figure}

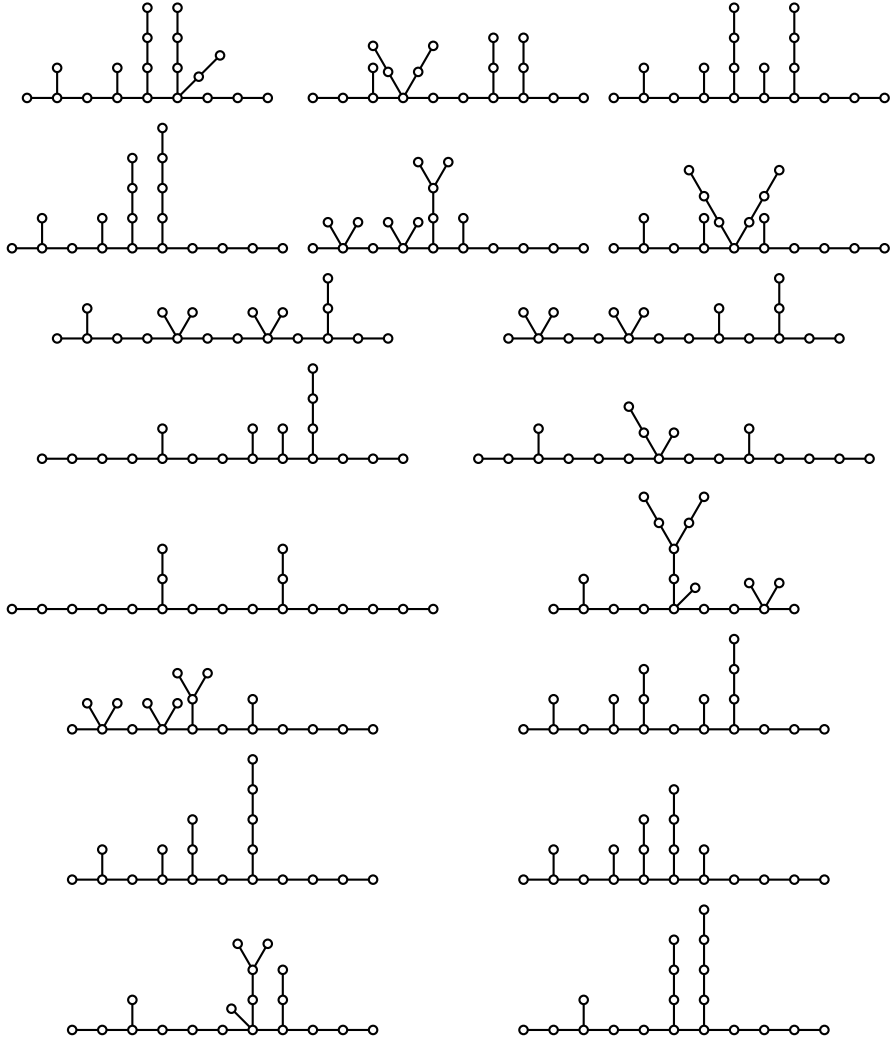
\begin{figure}[ht!]
\begin{center}
\begin{tikzpicture}[scale=0.4,style=thick]
\def\vr{2pt/0.5}

\coordinate(g19x1y1) at (0.5,-4);
\begin{scope}[shift=(g19x1y1)]
\foreach \i in {0,...,8}
{
	\coordinate(x\i) at (\i,0);
}

\draw(x0) -- (x1) -- (x2) -- (x3) -- (x4) -- (x5) -- (x6) -- (x7) -- (x8);

\coordinate(x9) at ({1+1*cos(90.0)}, {1*sin(90.0)});
\draw(x1) -- (x9);

\coordinate(x10) at ({3+1*cos(90.0)}, {1*sin(90.0)});
\draw(x3) -- (x10);

\coordinate(x11) at ({4+1*cos(90.0)}, {1*sin(90.0)});
\coordinate(x12) at ({4+2*cos(90.0)}, {2*sin(90.0)});
\coordinate(x13) at ({4+3*cos(90.0)}, {3*sin(90.0)});
\draw(x4) -- (x11) -- (x12) -- (x13);

\coordinate(x14) at ({5+1*cos(45.0)}, {1*sin(45.0)});
\coordinate(x15) at ({5+2*cos(45.0)}, {2*sin(45.0)});
\draw(x5) -- (x14) -- (x15);

\coordinate(x16) at ({5+1*cos(90.0)}, {1*sin(90.0)});
\coordinate(x17) at ({5+2*cos(90.0)}, {2*sin(90.0)});
\coordinate(x18) at ({5+3*cos(90.0)}, {3*sin(90.0)});
\draw(x5) -- (x16) -- (x17) -- (x18);

\foreach \i in {0,...,18}
{
	\draw (x\i)[fill=white] circle (\vr);
}
\end{scope}

\coordinate(g19x1y2) at (10.0,-4);
\begin{scope}[shift=(g19x1y2)]
\foreach \i in {0,...,9}
{
	\coordinate(x\i) at (\i,0);
}

\draw(x0) -- (x1) -- (x2) -- (x3) -- (x4) -- (x5) -- (x6) -- (x7) -- (x8) -- (x9);

\coordinate(x10) at ({2+1*cos(90.0)}, {1*sin(90.0)});
\draw(x2) -- (x10);

\coordinate(x11) at ({3+1*cos(60.0)}, {1*sin(60.0)});
\coordinate(x12) at ({3+2*cos(60.0)}, {2*sin(60.0)});
\draw(x3) -- (x11) -- (x12);

\coordinate(x13) at ({3+1*cos(120.0)}, {1*sin(120.0)});
\coordinate(x14) at ({3+2*cos(120.0)}, {2*sin(120.0)});
\draw(x3) -- (x13) -- (x14);

\coordinate(x15) at ({6+1*cos(90.0)}, {1*sin(90.0)});
\coordinate(x16) at ({6+2*cos(90.0)}, {2*sin(90.0)});
\draw(x6) -- (x15) -- (x16);

\coordinate(x17) at ({7+1*cos(90.0)}, {1*sin(90.0)});
\coordinate(x18) at ({7+2*cos(90.0)}, {2*sin(90.0)});
\draw(x7) -- (x17) -- (x18);

\foreach \i in {0,...,18}
{
	\draw (x\i)[fill=white] circle (\vr);
}
\end{scope}

\coordinate(g19x1y3) at (20.0,-4);
\begin{scope}[shift=(g19x1y3)]
\foreach \i in {0,...,9}
{
	\coordinate(x\i) at (\i,0);
}

\draw(x0) -- (x1) -- (x2) -- (x3) -- (x4) -- (x5) -- (x6) -- (x7) -- (x8) -- (x9);

\coordinate(x10) at ({1+1*cos(90.0)}, {1*sin(90.0)});
\draw(x1) -- (x10);

\coordinate(x11) at ({3+1*cos(90.0)}, {1*sin(90.0)});
\draw(x3) -- (x11);

\coordinate(x12) at ({4+1*cos(90.0)}, {1*sin(90.0)});
\coordinate(x13) at ({4+2*cos(90.0)}, {2*sin(90.0)});
\coordinate(x14) at ({4+3*cos(90.0)}, {3*sin(90.0)});
\draw(x4) -- (x12) -- (x13) -- (x14);

\coordinate(x15) at ({5+1*cos(90.0)}, {1*sin(90.0)});
\draw(x5) -- (x15);

\coordinate(x16) at ({6+1*cos(90.0)}, {1*sin(90.0)});
\coordinate(x17) at ({6+2*cos(90.0)}, {2*sin(90.0)});
\coordinate(x18) at ({6+3*cos(90.0)}, {3*sin(90.0)});
\draw(x6) -- (x16) -- (x17) -- (x18);

\foreach \i in {0,...,18}
{
	\draw (x\i)[fill=white] circle (\vr);
}
\end{scope}

\coordinate(g19x2y1) at (0.0,-9);
\begin{scope}[shift=(g19x2y1)]
\foreach \i in {0,...,9}
{
	\coordinate(x\i) at (\i,0);
}

\draw(x0) -- (x1) -- (x2) -- (x3) -- (x4) -- (x5) -- (x6) -- (x7) -- (x8) -- (x9);

\coordinate(x10) at ({1+1*cos(90.0)}, {1*sin(90.0)});
\draw(x1) -- (x10);

\coordinate(x11) at ({3+1*cos(90.0)}, {1*sin(90.0)});
\draw(x3) -- (x11);

\coordinate(x12) at ({4+1*cos(90.0)}, {1*sin(90.0)});
\coordinate(x13) at ({4+2*cos(90.0)}, {2*sin(90.0)});
\coordinate(x14) at ({4+3*cos(90.0)}, {3*sin(90.0)});
\draw(x4) -- (x12) -- (x13) -- (x14);

\coordinate(x15) at ({5+1*cos(90.0)}, {1*sin(90.0)});
\coordinate(x16) at ({5+2*cos(90.0)}, {2*sin(90.0)});
\coordinate(x17) at ({5+3*cos(90.0)}, {3*sin(90.0)});
\coordinate(x18) at ({5+4*cos(90.0)}, {4*sin(90.0)});
\draw(x5) -- (x15) -- (x16) -- (x17) -- (x18);

\foreach \i in {0,...,18}
{
	\draw (x\i)[fill=white] circle (\vr);
}
\end{scope}

\coordinate(g19x2y2) at (10.0,-9);
\begin{scope}[shift=(g19x2y2)]
\foreach \i in {0,...,9}
{
	\coordinate(x\i) at (\i,0);
}

\draw(x0) -- (x1) -- (x2) -- (x3) -- (x4) -- (x5) -- (x6) -- (x7) -- (x8) -- (x9);

\coordinate(x10) at ({1+1*cos(60.0)}, {1*sin(60.0)});
\draw(x1) -- (x10);

\coordinate(x11) at ({1+1*cos(120.0)}, {1*sin(120.0)});
\draw(x1) -- (x11);

\coordinate(x12) at ({3+1*cos(60.0)}, {1*sin(60.0)});
\draw(x3) -- (x12);

\coordinate(x13) at ({3+1*cos(120.0)}, {1*sin(120.0)});
\draw(x3) -- (x13);

\coordinate(x14) at ({4+1*cos(90.0)}, {1*sin(90.0)});
\draw(x4) -- (x14);

\coordinate(x15) at ({4+1*cos(90.0)}, {2*sin(90.0)});
\coordinate(x16) at ({4+1*cos(60.0)}, {2+sin(60.0)});
\coordinate(x17) at ({4+1*cos(120.0)}, {2+sin(120.0)});
\draw(x14) -- (x15) -- (x16);
\draw(x15) -- (x17);

\coordinate(x18) at ({5+1*cos(90.0)}, {1*sin(90.0)});
\draw(x5) -- (x18);

\foreach \i in {0,...,18}
{
	\draw (x\i)[fill=white] circle (\vr);
}
\end{scope}

\coordinate(g19x2y3) at (20.0,-9);
\begin{scope}[shift=(g19x2y3)]
\foreach \i in {0,...,9}
{
	\coordinate(x\i) at (\i,0);
}

\draw(x0) -- (x1) -- (x2) -- (x3) -- (x4) -- (x5) -- (x6) -- (x7) -- (x8) -- (x9);

\coordinate(x10) at ({1+1*cos(90.0)}, {1*sin(90.0)});
\draw(x1) -- (x10);

\coordinate(x11) at ({3+1*cos(90.0)}, {1*sin(90.0)});
\draw(x3) -- (x11);

\coordinate(x12) at ({4+1*cos(60.0)}, {1*sin(60.0)});
\coordinate(x13) at ({4+2*cos(60.0)}, {2*sin(60.0)});
\coordinate(x14) at ({4+3*cos(60.0)}, {3*sin(60.0)});
\draw(x4) -- (x12) -- (x13) -- (x14);

\coordinate(x15) at ({4+1*cos(120.0)}, {1*sin(120.0)});
\coordinate(x16) at ({4+2*cos(120.0)}, {2*sin(120.0)});
\coordinate(x17) at ({4+3*cos(120.0)}, {3*sin(120.0)});
\draw(x4) -- (x15) -- (x16) -- (x17);

\coordinate(x18) at ({5+1*cos(90.0)}, {1*sin(90.0)});
\draw(x5) -- (x18);

\foreach \i in {0,...,18}
{
	\draw (x\i)[fill=white] circle (\vr);
}
\end{scope}

\coordinate(g19x3y1) at (1.5,-12);
\begin{scope}[shift=(g19x3y1)]
\foreach \i in {0,...,11}
{
	\coordinate(x\i) at (\i,0);
}

\draw(x0) -- (x1) -- (x2) -- (x3) -- (x4) -- (x5) -- (x6) -- (x7) -- (x8) -- (x9) -- (x10) -- (x11);

\coordinate(x12) at ({1+1*cos(90.0)}, {1*sin(90.0)});
\draw(x1) -- (x12);

\coordinate(x13) at ({4+1*cos(60.0)}, {1*sin(60.0)});
\draw(x4) -- (x13);

\coordinate(x14) at ({4+1*cos(120.0)}, {1*sin(120.0)});
\draw(x4) -- (x14);

\coordinate(x15) at ({7+1*cos(60.0)}, {1*sin(60.0)});
\draw(x7) -- (x15);

\coordinate(x16) at ({7+1*cos(120.0)}, {1*sin(120.0)});
\draw(x7) -- (x16);

\coordinate(x17) at ({9+1*cos(90.0)}, {1*sin(90.0)});
\coordinate(x18) at ({9+2*cos(90.0)}, {2*sin(90.0)});
\draw(x9) -- (x17) -- (x18);

\foreach \i in {0,...,18}
{
	\draw (x\i)[fill=white] circle (\vr);
}
\end{scope}

\coordinate(g19x3y2) at (16.5,-12);
\begin{scope}[shift=(g19x3y2)]
\foreach \i in {0,...,11}
{
	\coordinate(x\i) at (\i,0);
}

\draw(x0) -- (x1) -- (x2) -- (x3) -- (x4) -- (x5) -- (x6) -- (x7) -- (x8) -- (x9) -- (x10) -- (x11);

\coordinate(x12) at ({1+1*cos(60.0)}, {1*sin(60.0)});
\draw(x1) -- (x12);

\coordinate(x13) at ({1+1*cos(120.0)}, {1*sin(120.0)});
\draw(x1) -- (x13);

\coordinate(x14) at ({4+1*cos(60.0)}, {1*sin(60.0)});
\draw(x4) -- (x14);

\coordinate(x15) at ({4+1*cos(120.0)}, {1*sin(120.0)});
\draw(x4) -- (x15);

\coordinate(x16) at ({7+1*cos(90.0)}, {1*sin(90.0)});
\draw(x7) -- (x16);

\coordinate(x17) at ({9+1*cos(90.0)}, {1*sin(90.0)});
\coordinate(x18) at ({9+2*cos(90.0)}, {2*sin(90.0)});
\draw(x9) -- (x17) -- (x18);

\foreach \i in {0,...,18}
{
	\draw (x\i)[fill=white] circle (\vr);
}
\end{scope}

\coordinate(g19x4y1) at (1.0,-16);
\begin{scope}[shift=(g19x4y1)]
\foreach \i in {0,...,12}
{
	\coordinate(x\i) at (\i,0);
}

\draw(x0) -- (x1) -- (x2) -- (x3) -- (x4) -- (x5) -- (x6) -- (x7) -- (x8) -- (x9) -- (x10) -- (x11) -- (x12);

\coordinate(x13) at ({4+1*cos(90.0)}, {1*sin(90.0)});
\draw(x4) -- (x13);

\coordinate(x14) at ({7+1*cos(90.0)}, {1*sin(90.0)});
\draw(x7) -- (x14);

\coordinate(x15) at ({8+1*cos(90.0)}, {1*sin(90.0)});
\draw(x8) -- (x15);

\coordinate(x16) at ({9+1*cos(90.0)}, {1*sin(90.0)});
\coordinate(x17) at ({9+2*cos(90.0)}, {2*sin(90.0)});
\coordinate(x18) at ({9+3*cos(90.0)}, {3*sin(90.0)});
\draw(x9) -- (x16) -- (x17) -- (x18);

\foreach \i in {0,...,18}
{
	\draw (x\i)[fill=white] circle (\vr);
}
\end{scope}

\coordinate(g19x4y2) at (15.5,-16);
\begin{scope}[shift=(g19x4y2)]
\foreach \i in {0,...,13}
{
	\coordinate(x\i) at (\i,0);
}

\draw(x0) -- (x1) -- (x2) -- (x3) -- (x4) -- (x5) -- (x6) -- (x7) -- (x8) -- (x9) -- (x10) -- (x11) -- (x12) -- (x13);

\coordinate(x14) at ({2+1*cos(90.0)}, {1*sin(90.0)});
\draw(x2) -- (x14);

\coordinate(x15) at ({6+1*cos(60.0)}, {1*sin(60.0)});
\draw(x6) -- (x15);

\coordinate(x16) at ({6+1*cos(120.0)}, {1*sin(120.0)});
\coordinate(x17) at ({6+2*cos(120.0)}, {2*sin(120.0)});
\draw(x6) -- (x16) -- (x17);

\coordinate(x18) at ({9+1*cos(90.0)}, {1*sin(90.0)});
\draw(x9) -- (x18);

\foreach \i in {0,...,18}
{
	\draw (x\i)[fill=white] circle (\vr);
}
\end{scope}

\coordinate(g19x5y1) at (0.0,-21);
\begin{scope}[shift=(g19x5y1)]
\foreach \i in {0,...,14}
{
	\coordinate(x\i) at (\i,0);
}

\draw(x0) -- (x1) -- (x2) -- (x3) -- (x4) -- (x5) -- (x6) -- (x7) -- (x8) -- (x9) -- (x10) -- (x11) -- (x12) -- (x13) -- (x14);

\coordinate(x15) at ({5+1*cos(90.0)}, {1*sin(90.0)});
\coordinate(x16) at ({5+2*cos(90.0)}, {2*sin(90.0)});
\draw(x5) -- (x15) -- (x16);

\coordinate(x17) at ({9+1*cos(90.0)}, {1*sin(90.0)});
\coordinate(x18) at ({9+2*cos(90.0)}, {2*sin(90.0)});
\draw(x9) -- (x17) -- (x18);

\foreach \i in {0,...,18}
{
	\draw (x\i)[fill=white] circle (\vr);
}
\end{scope}

\coordinate(g19x5y2) at (18.0,-21);
\begin{scope}[shift=(g19x5y2)]
\foreach \i in {0,...,8}
{
	\coordinate(x\i) at (\i,0);
}

\draw(x0) -- (x1) -- (x2) -- (x3) -- (x4) -- (x5) -- (x6) -- (x7) -- (x8);

\coordinate(x9) at ({1+1*cos(90.0)}, {1*sin(90.0)});
\draw(x1) -- (x9);

\coordinate(x10) at ({4+1*cos(45.0)}, {1*sin(45.0)});
\draw(x4) -- (x10);

\coordinate(x11) at ({4+1*cos(90.0)}, {1*sin(90.0)});
\coordinate(x12) at ({4+1*cos(90.0)}, {2*sin(90.0)});
\coordinate(x13) at ({4+1*cos(60.0)}, {2*sin(90.0)+ 1*sin(60.0)});
\coordinate(x14) at ({4+2*cos(60.0)}, {2*sin(90.0)+2*sin(60.0)});
\coordinate(x15) at ({4+1*cos(120.0)}, {2*sin(90.0)+1*sin(120.0)});
\coordinate(x16) at ({4+2*cos(120.0)}, {2*sin(90.0)+2*sin(120.0)});
\draw(x4) -- (x11) -- (x12) -- (x13) -- (x14);
\draw(x12) -- (x15) -- (x16);

\coordinate(x17) at ({7+1*cos(60.0)}, {1*sin(60.0)});
\draw(x7) -- (x17);

\coordinate(x18) at ({7+1*cos(120.0)}, {1*sin(120.0)});
\draw(x7) -- (x18);

\foreach \i in {0,...,18}
{
	\draw (x\i)[fill=white] circle (\vr);
}
\end{scope}

\coordinate(g19x6y1) at (2.0,-25);
\begin{scope}[shift=(g19x6y1)]
\foreach \i in {0,...,10}
{
	\coordinate(x\i) at (\i,0);
}

\draw(x0) -- (x1) -- (x2) -- (x3) -- (x4) -- (x5) -- (x6) -- (x7) -- (x8) -- (x9) -- (x10);

\coordinate(x11) at ({1+1*cos(60.0)}, {1*sin(60.0)});
\draw(x1) -- (x11);

\coordinate(x12) at ({1+1*cos(120.0)}, {1*sin(120.0)});
\draw(x1) -- (x12);

\coordinate(x13) at ({3+1*cos(60.0)}, {1*sin(60.0)});
\draw(x3) -- (x13);

\coordinate(x14) at ({3+1*cos(120.0)}, {1*sin(120.0)});
\draw(x3) -- (x14);

\coordinate(x15) at ({4+1*cos(90.0)}, {1*sin(90.0)});
\coordinate(x16) at ({4+1*cos(60.0)}, {1*sin(90.0)+1*sin(60.0)});
\coordinate(x17) at ({4+1*cos(120.0)}, {1*sin(90.0)+1*sin(120.0)});
\draw(x4) -- (x15) -- (x16);
\draw(x15) -- (x17);

\coordinate(x18) at ({6+1*cos(90.0)}, {1*sin(90.0)});
\draw(x6) -- (x18);

\foreach \i in {0,...,18}
{
	\draw (x\i)[fill=white] circle (\vr);
}
\end{scope}

\coordinate(g19x6y2) at (17.0,-25);
\begin{scope}[shift=(g19x6y2)]
\foreach \i in {0,...,10}
{
	\coordinate(x\i) at (\i,0);
}

\draw(x0) -- (x1) -- (x2) -- (x3) -- (x4) -- (x5) -- (x6) -- (x7) -- (x8) -- (x9) -- (x10);

\coordinate(x11) at ({1+1*cos(90.0)}, {1*sin(90.0)});
\draw(x1) -- (x11);

\coordinate(x12) at ({3+1*cos(90.0)}, {1*sin(90.0)});
\draw(x3) -- (x12);

\coordinate(x13) at ({4+1*cos(90.0)}, {1*sin(90.0)});
\coordinate(x14) at ({4+2*cos(90.0)}, {2*sin(90.0)});
\draw(x4) -- (x13) -- (x14);

\coordinate(x15) at ({6+1*cos(90.0)}, {1*sin(90.0)});
\draw(x6) -- (x15);

\coordinate(x16) at ({7+1*cos(90.0)}, {1*sin(90.0)});
\coordinate(x17) at ({7+2*cos(90.0)}, {2*sin(90.0)});
\coordinate(x18) at ({7+3*cos(90.0)}, {3*sin(90.0)});
\draw(x7) -- (x16) -- (x17) -- (x18);

\foreach \i in {0,...,18}
{
	\draw (x\i)[fill=white] circle (\vr);
}
\end{scope}

\coordinate(g19x7y1) at (2.0,-30);
\begin{scope}[shift=(g19x7y1)]
\foreach \i in {0,...,10}
{
	\coordinate(x\i) at (\i,0);
}

\draw(x0) -- (x1) -- (x2) -- (x3) -- (x4) -- (x5) -- (x6) -- (x7) -- (x8) -- (x9) -- (x10);

\coordinate(x11) at ({1+1*cos(90.0)}, {1*sin(90.0)});
\draw(x1) -- (x11);

\coordinate(x12) at ({3+1*cos(90.0)}, {1*sin(90.0)});
\draw(x3) -- (x12);

\coordinate(x13) at ({4+1*cos(90.0)}, {1*sin(90.0)});
\coordinate(x14) at ({4+2*cos(90.0)}, {2*sin(90.0)});
\draw(x4) -- (x13) -- (x14);

\coordinate(x15) at ({6+1*cos(90.0)}, {1*sin(90.0)});
\coordinate(x16) at ({6+2*cos(90.0)}, {2*sin(90.0)});
\coordinate(x17) at ({6+3*cos(90.0)}, {3*sin(90.0)});
\coordinate(x18) at ({6+4*cos(90.0)}, {4*sin(90.0)});
\draw(x6) -- (x15) -- (x16) -- (x17) -- (x18);

\foreach \i in {0,...,18}
{
	\draw (x\i)[fill=white] circle (\vr);
}
\end{scope}

\coordinate(g19x7y2) at (17.0,-30);
\begin{scope}[shift=(g19x7y2)]
\foreach \i in {0,...,10}
{
	\coordinate(x\i) at (\i,0);
}

\draw(x0) -- (x1) -- (x2) -- (x3) -- (x4) -- (x5) -- (x6) -- (x7) -- (x8) -- (x9) -- (x10);

\coordinate(x11) at ({1+1*cos(90.0)}, {1*sin(90.0)});
\draw(x1) -- (x11);

\coordinate(x12) at ({3+1*cos(90.0)}, {1*sin(90.0)});
\draw(x3) -- (x12);

\coordinate(x13) at ({4+1*cos(90.0)}, {1*sin(90.0)});
\coordinate(x14) at ({4+2*cos(90.0)}, {2*sin(90.0)});
\draw(x4) -- (x13) -- (x14);

\coordinate(x15) at ({5+1*cos(90.0)}, {1*sin(90.0)});
\coordinate(x16) at ({5+2*cos(90.0)}, {2*sin(90.0)});
\coordinate(x17) at ({5+3*cos(90.0)}, {3*sin(90.0)});
\draw(x5) -- (x15) -- (x16) -- (x17);

\coordinate(x18) at ({6+1*cos(90.0)}, {1*sin(90.0)});
\draw(x6) -- (x18);

\foreach \i in {0,...,18}
{
	\draw (x\i)[fill=white] circle (\vr);
}
\end{scope}

\coordinate(g19x8y1) at (2.0,-35);
\begin{scope}[shift=(g19x8y1)]
\foreach \i in {0,...,10}
{
	\coordinate(x\i) at (\i,0);
}

\draw(x0) -- (x1) -- (x2) -- (x3) -- (x4) -- (x5) -- (x6) -- (x7) -- (x8) -- (x9) -- (x10);

\coordinate(x11) at ({2+1*cos(90.0)}, {1*sin(90.0)});
\draw(x2) -- (x11);

\coordinate(x12) at ({6+1*cos(135.0)}, {1*sin(135.0)});
\draw(x6) -- (x12);

\coordinate(x13) at ({6+1*cos(90.0)}, {1*sin(90.0)});
\coordinate(x14) at ({6+2*cos(90.0)}, {2*sin(90.0)});
\coordinate(x15) at ({6+1*cos(60.0)}, {2*sin(90.0)+1*sin(60.0)});
\coordinate(x16) at ({6+1*cos(120.0)}, {2*sin(90.0)+1*sin(120.0)});
\draw(x6) -- (x13) -- (x14) -- (x15);
\draw(x14) -- (x16);

\coordinate(x17) at ({7+1*cos(90.0)}, {1*sin(90.0)});
\coordinate(x18) at ({7+2*cos(90.0)}, {2*sin(90.0)});
\draw(x7) -- (x17) -- (x18);

\foreach \i in {0,...,18}
{
	\draw (x\i)[fill=white] circle (\vr);
}
\end{scope}

\coordinate(g19x8y2) at (17.0,-35);
\begin{scope}[shift=(g19x8y2)]
\foreach \i in {0,...,10}
{
	\coordinate(x\i) at (\i,0);
}

\draw(x0) -- (x1) -- (x2) -- (x3) -- (x4) -- (x5) -- (x6) -- (x7) -- (x8) -- (x9) -- (x10);

\coordinate(x11) at ({2+1*cos(90.0)}, {1*sin(90.0)});
\draw(x2) -- (x11);

\coordinate(x12) at ({5+1*cos(90.0)}, {1*sin(90.0)});
\coordinate(x13) at ({5+2*cos(90.0)}, {2*sin(90.0)});
\coordinate(x14) at ({5+3*cos(90.0)}, {3*sin(90.0)});
\draw(x5) -- (x12) -- (x13) -- (x14);

\coordinate(x15) at ({6+1*cos(90.0)}, {1*sin(90.0)});
\coordinate(x16) at ({6+2*cos(90.0)}, {2*sin(90.0)});
\coordinate(x17) at ({6+3*cos(90.0)}, {3*sin(90.0)});
\coordinate(x18) at ({6+4*cos(90.0)}, {4*sin(90.0)});
\draw(x6) -- (x15) -- (x16) -- (x17) -- (x18);

\foreach \i in {0,...,18}
{
	\draw (x\i)[fill=white] circle (\vr);
}
\end{scope}

\end{tikzpicture}
\end{center}
\caption{All critical trees on 19 vertices}
\label{fig:critical-trees-19}
\end{figure}

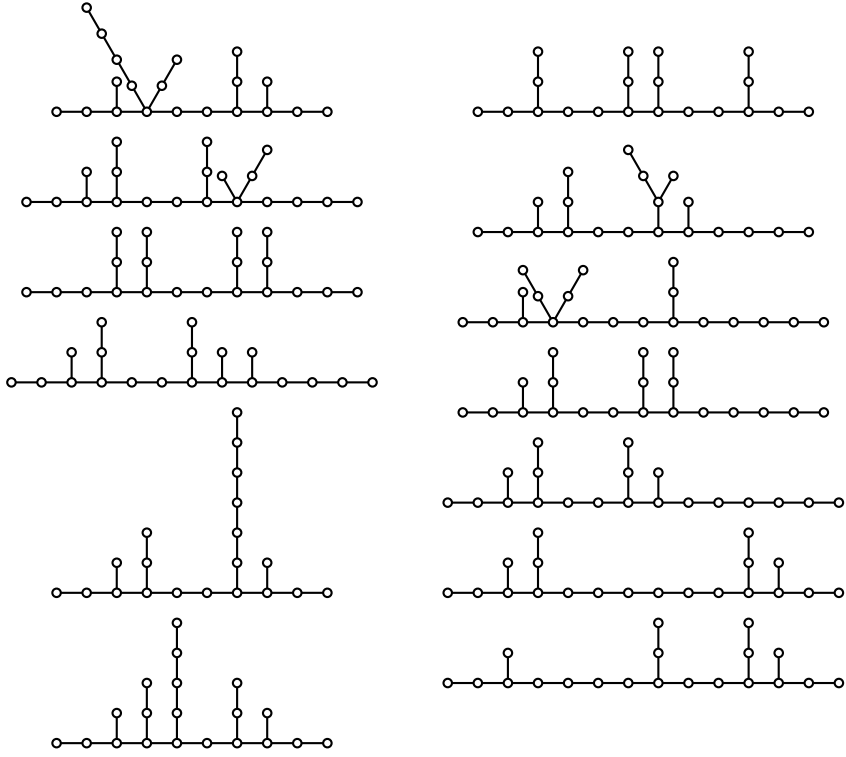
\begin{figure}[ht!]
\begin{center}
\begin{tikzpicture}[scale=0.4,style=thick]
\def\vr{2pt/0.5}

\coordinate(g20x1y1) at (2.5,-5);
\begin{scope}[shift=(g20x1y1)]
\foreach \i in {0,...,9}
{
	\coordinate(x\i) at (\i,0);
}

\draw(x0) -- (x1) -- (x2) -- (x3) -- (x4) -- (x5) -- (x6) -- (x7) -- (x8) -- (x9);

\coordinate(x10) at ({2+1*cos(90.0)}, {1*sin(90.0)});
\draw(x2) -- (x10);

\coordinate(x11) at ({3+1*cos(60.0)}, {1*sin(60.0)});
\coordinate(x12) at ({3+2*cos(60.0)}, {2*sin(60.0)});
\draw(x3) -- (x11) -- (x12);

\coordinate(x13) at ({3+1*cos(120.0)}, {1*sin(120.0)});
\coordinate(x14) at ({3+2*cos(120.0)}, {2*sin(120.0)});
\coordinate(x15) at ({3+3*cos(120.0)}, {3*sin(120.0)});
\coordinate(x16) at ({3+4*cos(120.0)}, {4*sin(120.0)});
\draw(x3) -- (x13) -- (x14) -- (x15) -- (x16);

\coordinate(x17) at ({6+1*cos(90.0)}, {1*sin(90.0)});
\coordinate(x18) at ({6+2*cos(90.0)}, {2*sin(90.0)});
\draw(x6) -- (x17) -- (x18);

\coordinate(x19) at ({7+1*cos(90.0)}, {1*sin(90.0)});
\draw(x7) -- (x19);

\foreach \i in {0,...,19}
{
	\draw (x\i)[fill=white] circle (\vr);
}
\end{scope}

\coordinate(g20x1y2) at (16.5,-5);
\begin{scope}[shift=(g20x1y2)]
\foreach \i in {0,...,11}
{
	\coordinate(x\i) at (\i,0);
}

\draw(x0) -- (x1) -- (x2) -- (x3) -- (x4) -- (x5) -- (x6) -- (x7) -- (x8) -- (x9) -- (x10) -- (x11);

\coordinate(x12) at ({2+1*cos(90.0)}, {1*sin(90.0)});
\coordinate(x13) at ({2+2*cos(90.0)}, {2*sin(90.0)});
\draw(x2) -- (x12) -- (x13);

\coordinate(x14) at ({5+1*cos(90.0)}, {1*sin(90.0)});
\coordinate(x15) at ({5+2*cos(90.0)}, {2*sin(90.0)});
\draw(x5) -- (x14) -- (x15);

\coordinate(x16) at ({6+1*cos(90.0)}, {1*sin(90.0)});
\coordinate(x17) at ({6+2*cos(90.0)}, {2*sin(90.0)});
\draw(x6) -- (x16) -- (x17);

\coordinate(x18) at ({9+1*cos(90.0)}, {1*sin(90.0)});
\coordinate(x19) at ({9+2*cos(90.0)}, {2*sin(90.0)});
\draw(x9) -- (x18) -- (x19);

\foreach \i in {0,...,19}
{
	\draw (x\i)[fill=white] circle (\vr);
}
\end{scope}

\coordinate(g20x2y1) at (1.5,-8);
\begin{scope}[shift=(g20x2y1)]
\foreach \i in {0,...,11}
{
	\coordinate(x\i) at (\i,0);
}

\draw(x0) -- (x1) -- (x2) -- (x3) -- (x4) -- (x5) -- (x6) -- (x7) -- (x8) -- (x9) -- (x10) -- (x11);

\coordinate(x12) at ({2+1*cos(90.0)}, {1*sin(90.0)});
\draw(x2) -- (x12);

\coordinate(x13) at ({3+1*cos(90.0)}, {1*sin(90.0)});
\coordinate(x14) at ({3+2*cos(90.0)}, {2*sin(90.0)});
\draw(x3) -- (x13) -- (x14);

\coordinate(x15) at ({6+1*cos(90.0)}, {1*sin(90.0)});
\coordinate(x16) at ({6+2*cos(90.0)}, {2*sin(90.0)});
\draw(x6) -- (x15) -- (x16);

\coordinate(x17) at ({7+1*cos(60.0)}, {1*sin(60.0)});
\coordinate(x18) at ({7+2*cos(60.0)}, {2*sin(60.0)});
\draw(x7) -- (x17) -- (x18);

\coordinate(x19) at ({7+1*cos(120.0)}, {1*sin(120.0)});
\draw(x7) -- (x19);

\foreach \i in {0,...,19}
{
	\draw (x\i)[fill=white] circle (\vr);
}
\end{scope}

\coordinate(g20x2y2) at (16.5,-9);
\begin{scope}[shift=(g20x2y2)]
\foreach \i in {0,...,11}
{
	\coordinate(x\i) at (\i,0);
}

\draw(x0) -- (x1) -- (x2) -- (x3) -- (x4) -- (x5) -- (x6) -- (x7) -- (x8) -- (x9) -- (x10) -- (x11);

\coordinate(x12) at ({2+1*cos(90.0)}, {1*sin(90.0)});
\draw(x2) -- (x12);

\coordinate(x13) at ({3+1*cos(90.0)}, {1*sin(90.0)});
\coordinate(x14) at ({3+2*cos(90.0)}, {2*sin(90.0)});
\draw(x3) -- (x13) -- (x14);

\coordinate(x15) at ({6+1*cos(90.0)}, {1*sin(90.0)});
\coordinate(x16) at ({6+1*cos(90.0)+cos(120)}, {sin(90.0)+sin(120)});
\coordinate(x17) at ({6+3*cos(90.0)+2*cos(120)}, {sin(90.0)+2*sin(120)});
\coordinate(x18) at ({6+4*cos(90.0)+cos(60)}, {sin(90.0)+sin(60)});
\draw(x6) -- (x15) -- (x16) -- (x17);
\draw(x15) -- (x18);

\coordinate(x19) at ({7+1*cos(90.0)}, {1*sin(90.0)});
\draw(x7) -- (x19);

\foreach \i in {0,...,19}
{
	\draw (x\i)[fill=white] circle (\vr);
}
\end{scope}

\coordinate(g20x3y1) at (1.5,-11);
\begin{scope}[shift=(g20x3y1)]
\foreach \i in {0,...,11}
{
	\coordinate(x\i) at (\i,0);
}

\draw(x0) -- (x1) -- (x2) -- (x3) -- (x4) -- (x5) -- (x6) -- (x7) -- (x8) -- (x9) -- (x10) -- (x11);

\coordinate(x12) at ({3+1*cos(90.0)}, {1*sin(90.0)});
\coordinate(x13) at ({3+2*cos(90.0)}, {2*sin(90.0)});
\draw(x3) -- (x12) -- (x13);

\coordinate(x14) at ({4+1*cos(90.0)}, {1*sin(90.0)});
\coordinate(x15) at ({4+2*cos(90.0)}, {2*sin(90.0)});
\draw(x4) -- (x14) -- (x15);

\coordinate(x16) at ({7+1*cos(90.0)}, {1*sin(90.0)});
\coordinate(x17) at ({7+2*cos(90.0)}, {2*sin(90.0)});
\draw(x7) -- (x16) -- (x17);

\coordinate(x18) at ({8+1*cos(90.0)}, {1*sin(90.0)});
\coordinate(x19) at ({8+2*cos(90.0)}, {2*sin(90.0)});
\draw(x8) -- (x18) -- (x19);

\foreach \i in {0,...,19}
{
	\draw (x\i)[fill=white] circle (\vr);
}
\end{scope}

\coordinate(g20x3y2) at (16.0,-12);
\begin{scope}[shift=(g20x3y2)]
\foreach \i in {0,...,12}
{
	\coordinate(x\i) at (\i,0);
}

\draw(x0) -- (x1) -- (x2) -- (x3) -- (x4) -- (x5) -- (x6) -- (x7) -- (x8) -- (x9) -- (x10) -- (x11) -- (x12);

\coordinate(x13) at ({2+1*cos(90.0)}, {1*sin(90.0)});
\draw(x2) -- (x13);

\coordinate(x14) at ({3+1*cos(60.0)}, {1*sin(60.0)});
\coordinate(x15) at ({3+2*cos(60.0)}, {2*sin(60.0)});
\draw(x3) -- (x14) -- (x15);

\coordinate(x16) at ({3+1*cos(120.0)}, {1*sin(120.0)});
\coordinate(x17) at ({3+2*cos(120.0)}, {2*sin(120.0)});
\draw(x3) -- (x16) -- (x17);

\coordinate(x18) at ({7+1*cos(90.0)}, {1*sin(90.0)});
\coordinate(x19) at ({7+2*cos(90.0)}, {2*sin(90.0)});
\draw(x7) -- (x18) -- (x19);

\foreach \i in {0,...,19}
{
	\draw (x\i)[fill=white] circle (\vr);
}
\end{scope}

\coordinate(g20x4y1) at (1.0,-14);
\begin{scope}[shift=(g20x4y1)]
\foreach \i in {0,...,12}
{
	\coordinate(x\i) at (\i,0);
}

\draw(x0) -- (x1) -- (x2) -- (x3) -- (x4) -- (x5) -- (x6) -- (x7) -- (x8) -- (x9) -- (x10) -- (x11) -- (x12);

\coordinate(x13) at ({2+1*cos(90.0)}, {1*sin(90.0)});
\draw(x2) -- (x13);

\coordinate(x14) at ({3+1*cos(90.0)}, {1*sin(90.0)});
\coordinate(x15) at ({3+2*cos(90.0)}, {2*sin(90.0)});
\draw(x3) -- (x14) -- (x15);

\coordinate(x16) at ({6+1*cos(90.0)}, {1*sin(90.0)});
\coordinate(x17) at ({6+2*cos(90.0)}, {2*sin(90.0)});
\draw(x6) -- (x16) -- (x17);

\coordinate(x18) at ({7+1*cos(90.0)}, {1*sin(90.0)});
\draw(x7) -- (x18);

\coordinate(x19) at ({8+1*cos(90.0)}, {1*sin(90.0)});
\draw(x8) -- (x19);

\foreach \i in {0,...,19}
{
	\draw (x\i)[fill=white] circle (\vr);
}
\end{scope}

\coordinate(g20x4y2) at (16.0,-15);
\begin{scope}[shift=(g20x4y2)]
\foreach \i in {0,...,12}
{
	\coordinate(x\i) at (\i,0);
}

\draw(x0) -- (x1) -- (x2) -- (x3) -- (x4) -- (x5) -- (x6) -- (x7) -- (x8) -- (x9) -- (x10) -- (x11) -- (x12);

\coordinate(x13) at ({2+1*cos(90.0)}, {1*sin(90.0)});
\draw(x2) -- (x13);

\coordinate(x14) at ({3+1*cos(90.0)}, {1*sin(90.0)});
\coordinate(x15) at ({3+2*cos(90.0)}, {2*sin(90.0)});
\draw(x3) -- (x14) -- (x15);

\coordinate(x16) at ({6+1*cos(90.0)}, {1*sin(90.0)});
\coordinate(x17) at ({6+2*cos(90.0)}, {2*sin(90.0)});
\draw(x6) -- (x16) -- (x17);

\coordinate(x18) at ({7+1*cos(90.0)}, {1*sin(90.0)});
\coordinate(x19) at ({7+2*cos(90.0)}, {2*sin(90.0)});
\draw(x7) -- (x18) -- (x19);

\foreach \i in {0,...,19}
{
	\draw (x\i)[fill=white] circle (\vr);
}
\end{scope}

\coordinate(g20x5y1) at (2.5,-21);
\begin{scope}[shift=(g20x5y1)]
\foreach \i in {0,...,9}
{
	\coordinate(x\i) at (\i,0);
}

\draw(x0) -- (x1) -- (x2) -- (x3) -- (x4) -- (x5) -- (x6) -- (x7) -- (x8) -- (x9);

\coordinate(x10) at ({2+1*cos(90.0)}, {1*sin(90.0)});
\draw(x2) -- (x10);

\coordinate(x11) at ({3+1*cos(90.0)}, {1*sin(90.0)});
\coordinate(x12) at ({3+2*cos(90.0)}, {2*sin(90.0)});
\draw(x3) -- (x11) -- (x12);

\coordinate(x13) at ({6+1*cos(90.0)}, {1*sin(90.0)});
\coordinate(x14) at ({6+2*cos(90.0)}, {2*sin(90.0)});
\coordinate(x15) at ({6+3*cos(90.0)}, {3*sin(90.0)});
\coordinate(x16) at ({6+4*cos(90.0)}, {4*sin(90.0)});
\coordinate(x17) at ({6+5*cos(90.0)}, {5*sin(90.0)});
\coordinate(x18) at ({6+6*cos(90.0)}, {6*sin(90.0)});
\draw(x6) -- (x13) -- (x14) -- (x15) -- (x16) -- (x17) -- (x18);

\coordinate(x19) at ({7+1*cos(90.0)}, {1*sin(90.0)});
\draw(x7) -- (x19);

\foreach \i in {0,...,19}
{
	\draw (x\i)[fill=white] circle (\vr);
}
\end{scope}

\coordinate(g20x5y2) at (15.5,-21);
\begin{scope}[shift=(g20x5y2)]
\foreach \i in {0,...,13}
{
	\coordinate(x\i) at (\i,0);
}

\draw(x0) -- (x1) -- (x2) -- (x3) -- (x4) -- (x5) -- (x6) -- (x7) -- (x8) -- (x9) -- (x10) -- (x11) -- (x12) -- (x13);

\coordinate(x14) at ({2+1*cos(90.0)}, {1*sin(90.0)});
\draw(x2) -- (x14);

\coordinate(x15) at ({3+1*cos(90.0)}, {1*sin(90.0)});
\coordinate(x16) at ({3+2*cos(90.0)}, {2*sin(90.0)});
\draw(x3) -- (x15) -- (x16);

\coordinate(x17) at ({10+1*cos(90.0)}, {1*sin(90.0)});
\coordinate(x18) at ({10+2*cos(90.0)}, {2*sin(90.0)});
\draw(x10) -- (x17) -- (x18);

\coordinate(x19) at ({11+1*cos(90.0)}, {1*sin(90.0)});
\draw(x11) -- (x19);

\foreach \i in {0,...,19}
{
	\draw (x\i)[fill=white] circle (\vr);
}
\end{scope}

\coordinate(g20x6y1) at (15.5,-24);
\begin{scope}[shift=(g20x6y1)]
\foreach \i in {0,...,13}
{
	\coordinate(x\i) at (\i,0);
}

\draw(x0) -- (x1) -- (x2) -- (x3) -- (x4) -- (x5) -- (x6) -- (x7) -- (x8) -- (x9) -- (x10) -- (x11) -- (x12) -- (x13);

\coordinate(x14) at ({2+1*cos(90.0)}, {1*sin(90.0)});
\draw(x2) -- (x14);

\coordinate(x15) at ({7+1*cos(90.0)}, {1*sin(90.0)});
\coordinate(x16) at ({7+2*cos(90.0)}, {2*sin(90.0)});
\draw(x7) -- (x15) -- (x16);

\coordinate(x17) at ({10+1*cos(90.0)}, {1*sin(90.0)});
\coordinate(x18) at ({10+2*cos(90.0)}, {2*sin(90.0)});
\draw(x10) -- (x17) -- (x18);

\coordinate(x19) at ({11+1*cos(90.0)}, {1*sin(90.0)});
\draw(x11) -- (x19);

\foreach \i in {0,...,19}
{
	\draw (x\i)[fill=white] circle (\vr);
}
\end{scope}

\coordinate(g20x6y2) at (15.5,-18);
\begin{scope}[shift=(g20x6y2)]
\foreach \i in {0,...,13}
{
	\coordinate(x\i) at (\i,0);
}

\draw(x0) -- (x1) -- (x2) -- (x3) -- (x4) -- (x5) -- (x6) -- (x7) -- (x8) -- (x9) -- (x10) -- (x11) -- (x12) -- (x13);

\coordinate(x14) at ({2+1*cos(90.0)}, {1*sin(90.0)});
\draw(x2) -- (x14);

\coordinate(x15) at ({3+1*cos(90.0)}, {1*sin(90.0)});
\coordinate(x16) at ({3+2*cos(90.0)}, {2*sin(90.0)});
\draw(x3) -- (x15) -- (x16);

\coordinate(x17) at ({6+1*cos(90.0)}, {1*sin(90.0)});
\coordinate(x18) at ({6+2*cos(90.0)}, {2*sin(90.0)});
\draw(x6) -- (x17) -- (x18);

\coordinate(x19) at ({7+1*cos(90.0)}, {1*sin(90.0)});
\draw(x7) -- (x19);

\foreach \i in {0,...,19}
{
	\draw (x\i)[fill=white] circle (\vr);
}
\end{scope}

\coordinate(g20x7y1) at (2.5,-26);
\begin{scope}[shift=(g20x7y1)]
\foreach \i in {0,...,9}
{
	\coordinate(x\i) at (\i,0);
}

\draw(x0) -- (x1) -- (x2) -- (x3) -- (x4) -- (x5) -- (x6) -- (x7) -- (x8) -- (x9);

\coordinate(x10) at ({2+1*cos(90.0)}, {1*sin(90.0)});
\draw(x2) -- (x10);

\coordinate(x11) at ({3+1*cos(90.0)}, {1*sin(90.0)});
\coordinate(x12) at ({3+2*cos(90.0)}, {2*sin(90.0)});
\draw(x3) -- (x11) -- (x12);

\coordinate(x13) at ({4+1*cos(90.0)}, {1*sin(90.0)});
\coordinate(x14) at ({4+2*cos(90.0)}, {2*sin(90.0)});
\coordinate(x15) at ({4+3*cos(90.0)}, {3*sin(90.0)});
\coordinate(x16) at ({4+4*cos(90.0)}, {4*sin(90.0)});
\draw(x4) -- (x13) -- (x14) -- (x15) -- (x16);

\coordinate(x17) at ({6+1*cos(90.0)}, {1*sin(90.0)});
\coordinate(x18) at ({6+2*cos(90.0)}, {2*sin(90.0)});
\draw(x6) -- (x17) -- (x18);

\coordinate(x19) at ({7+1*cos(90.0)}, {1*sin(90.0)});
\draw(x7) -- (x19);

\foreach \i in {0,...,19}
{
	\draw (x\i)[fill=white] circle (\vr);
}
\end{scope}

\end{tikzpicture}
\end{center}
\caption{All critical trees on 20 vertices}
\label{fig:critical-trees-20}
\end{figure}

In~\cite{bujtas-2015c} it was reported that there are no $\gamma_g$-critical trees on up to 12 vertices, there are two $\gamma_g$-critical trees on 13 vertices, there are no such trees on 14 and 15 vertices, there is only one such tree on 16 vertices, but there are ten on 17 vertices. So the variety of $\gamma_g$-critical trees appeared to be quite small. We  were now able to extend our computations and obtain all $\gamma_g$-critical trees on up to 20 vertices. It turned out that the variety is somehow larger than expected! The obtained trees are drawn in Figs.~\ref{fig:critical-trees-18}, \ref{fig:critical-trees-19}, and~\ref{fig:critical-trees-20}.

From the list of $\gamma_g$-critical trees on up to 20 vertices we were not able to formulate a conjecture asserting a characterization of $\gamma_g$-critical trees. This is in contrast with the fact that in~\cite{henning-2017b} a conjecture asserting the structure of the $3/5$-trees has been posed.

\section{Extended Cutting Lemma}
\label{sec:extended-cut-lemma}

In this section we extend the Cutting Lemma so that it involves also the Staller-pass game. We begin with the following:

\begin{lemma}
\label{lem:Cutting Lemma II}
Let $G$ be a graph, and let $B,C \subseteq V(G)$ where $C \subseteq B$. If $uv$ is an edge of $G$ and $\{u,v\} \subseteq C$, then
\[
\dstart(G_{uv}|B) \le \dstart(G|C) \hspace*{0.5cm} \mbox{and} \hspace*{0.5cm} \sstart(G_{uv}|B) \le \sstart(G|C).
\]
\end{lemma}

\proof Since the vertex $u'$ and its neighbor $v$ in $G_{uv}$ are both dominated in $G_{uv}|B$, the vertex $u'$ plays no further role in the game and can be deleted. Analogously, the vertex $v'$ plays no further role in the game played on $G_{uv}|B$ and can be deleted. Thus, $\dstart(G_{uv}|B) = \dstart(G - uv|B)$. Moreover, since both $u$ and $v$ are dominated in $G|C$, the edge $uv$ plays no further role in the game and can be deleted. Thus, $\dstart(G|C) = \dstart(G-uv|C)$. Since $C \subseteq B$, we can apply the Continuation Principle to the partially dominated graphs $G - uv|B$ and $G-uv|C$, yielding
\[
\dstart(G_{uv}|B) = \dstart(G - uv|B) \le \dstart(G-uv|C) = \dstart(G|C).
\]
Analogously, $\sstart(G_{uv}|B) \le \sstart(G|C)$.
\qed

For the next lemma, we recall that the Staller-pass game is the domination game in which \St is allowed, but not obligated, to skip exactly one move in the course of the game. Further, we recall that $\dstartpass(G)$ (resp. $\sstartpass(G)$) is the size of the dominating set produced under optimal play when Dominator (resp. Staller) starts the Staller-pass game. The turn when \St passes does not count as a move.

We shall need the following property of subsets of vertices of a graph.

\begin{unnumbered}{Definition (Subset Property)}
Let $G$ be a graph, and let $uv$ be an edge of $G$. For two subsets of vertices $B$ and $C$, we say that the ordered pair $(B,C)$ satisfies the \Buvprop\ if either $C \subseteq B$ or $\{u,v\} \subseteq C$ and $C \setminus B = \{w\}$ where $w \in \{u,v\}$.
\end{unnumbered}

We are now in a position to prove the following.

\begin{lemma}
\label{lem:Cutting Lemma III}
Let $G$ be a graph and let $uv$ be an edge of $G$. If
$B,C \subseteq V(G)$ and $(B,C)$ satisfies the \uvprop, then
\[
\dstart(G_{uv}|B) \le 1 + \dstartpass(G|C) \hspace*{0.5cm} \mbox{and} \hspace*{0.5cm} \sstart(G_{uv}|B) \le 1 + \sstartpass(G|C).
\]
\end{lemma}

\proof
We proceed by induction on $|V(G_{uv}) \setminus B|$. If $|V(G_{uv}) \setminus B| = 0$, then $\dstart(G_{uv}|B) = 0$ and $\dstartpass(G|C) \ge 0$. Further,  $\sstart(G_{uv}|B) = 0$ and $\sstartpass(G|C) \ge 0$. Thus the desired result holds. This establishes the base case. Suppose that $|V(G_{uv}) \setminus B| \ge 1$, and so the game is not yet complete in $G_{uv}|B$.
We prove first the upper bound on the S-game.

\begin{unnumbered}{Claim~A}
$\sstart(G_{uv}|B) \le 1 + \sstartpass(G|C)$.
\end{unnumbered}
\proof
We consider two cases. Suppose firstly that \emph{Staller has an optimal move $x$ in $G_{uv}|B$ which when played dominates a new vertex that does not belong to $C$.} Let $B' = B \cup N_{G_{uv}}[x]$ and $C' = C \cup N_G[x]$. If $x \notin \{u,v\}$, then $N_{G_{uv}}[x] = N_G[x]$ and since the pair $(B,C)$ satisfies the \uvprop, the pair $(B',C')$ also does.
If $x = u$, then $N_G[x] \setminus N_{G_{uv}}[x] \subseteq \{v\}$. Then $C'\setminus B'\subseteq \{v\}\cup(C\setminus B)$ and since $(B,C)$ satisfies the \uvprop, $C'\setminus B'\subseteq \{u,v\}$. But as $B'$ now contains $u$, we note that $C'\setminus B'\subseteq \{v\}$. Both $u$ and $v$ are in $N_G[u]\subseteq C'$, so the pair $(B',C')$ satisfies the \uvprop.
A symmetric argument can be made if $x=v$. Thus, by induction, we note that
\[
\begin{array}{lcl}
\sstart(G_{uv}|B) & \stackrel{(x {\rm \, optimal \, move})}{=} & 1 + \dstart(G_{uv}|B \cup N_{G_{uv}}[x]) \1 \\
& \stackrel{({\rm by \, induction})}{\le} & 1 + \left( 1 + \dstartpass(G|C \cup N_G[x]) \right) \2 \\
& \le & 1 + \sstartpass(G|C),
\end{array}
\]
where the last inequality follows from the observation that the vertex $x$ may not be an optimal move for \St in the Staller-pass game played in $G|C$.

Suppose secondly that \emph{Staller's optimal move in $G_{uv}|B$ is a vertex $x$ which when played dominates only one new vertex $w$ and this vertex belongs to $C \setminus B$.} In this case, we note that $\{u,v\} \subseteq C$ and $C \setminus B = \{w\}$ where $w \in \{u,v\}$. In particular, $C \subseteq B \cup \{w\}$. The following therefore holds.
\[
\begin{array}{lcl}
\sstart(G_{uv}|B) & \stackrel{(x {\rm \, optimal \, move})}{=} & 1 + \dstart(G_{uv}|B \cup \{w\}) \1 \\
& \stackrel{({\rm Lemma~\ref{lem:Cutting Lemma II})}}{\le} & 1 + \dstart(G|C)  \2 \\
& \le & 1 + \sstartpass(G|C),
\end{array}
\]
where the last inequality follows from the observation that passing may not be an optimal move for \St in the Staller-pass game played in $G|C$. This completes the proof of Claim~A.
\smallqed

\smallskip
We prove next the upper bound on the Dominator-start game.

\begin{unnumbered}{Claim~B}
$\dstart(G_{uv}|B) \le 1 + \dstartpass(G|C)$.
\end{unnumbered}

\proof
Let $x$ be an optimal move for \D in the Staller-pass game played in $G|C$. Let $B' = B \cup N_{G_{uv}}[x]$ and $C' = C \cup N_G[x]$. For $x \in \{u,v\}$, we let $\xbar =  \{u,v\} \setminus \{x\}$. Thus if $x = u$, then $\xbar = v$, while if $x = v$, then $\xbar = u$. We consider three cases.

\medskip
\emph{Case~1. $x \notin \{u,v\}$.} In this case, $N_{G_{uv}}[x] = N_G[x]$. Since the pair $(B,C)$ satisfies the \uvprop, the pair $(B',C')$ also satisfies the \uvprop. If $x$ is a legal move in $G_{uv}|B$, then by induction, we note that
\[
\begin{array}{lcl}
1 + \dstartpass(G|C) & \stackrel{(x {\rm \, optimal \, move})}{=} & 1 + (1 + \sstartpass(G|C \cup N_G[x])) \2 \\
& \stackrel{({\rm by \, induction})}{\ge} & 1 + \sstart(G_{uv}|B \cup N_{G_{uv}}[x]) \2 \\
& \ge & \dstart(G_{uv}|B),
\end{array}
\]

\noindent
where the last inequality follows from the observation that the vertex $x$ may not be an optimal move for \D in the game played in $G_{uv}|B$. Assume now that $x$ is not a legal move in $G_{uv}|B$. This implies that $B \cup N_{G_{uv}}[x] = B$. In this case, \D plays any move $w$ in $G_{uv}|B$ that is legal. The following now holds.
\[
\begin{array}{lcl}
1 + \dstartpass(G|C) & \stackrel{(x {\rm \, optimal \, move})}{=} & 1 + (1 + \sstartpass(G|C \cup N_G[x])) \1 \\
& \stackrel{({\rm by \, induction})}{\ge} & 1 + \sstart(G_{uv}|B \cup N_{G_{uv}}[x] \cup N_{G_{uv}}[w]) \2 \\
& = & 1 + \sstart(G_{uv}|B \cup N_{G_{uv}}[w]) \2 \\
& \ge & \dstart(G_{uv}|B),
\end{array}
\]
\noindent
where the last inequality follows from the observation that the vertex $w$ may not be an optimal move for \D in the game played in $G_{uv}|B$. This completes Case~1.

\medskip
\emph{Case~2. $x \in \{u,v\}$ and $\xbar \in B$.} Thus, either $x = u$ and $v \in B$ or $x = v$ and $u \in B$. Renaming $u$ and $v$ if necessary, we may assume that $x = u$ and $v \in B$. Hence, $C \cup N_G[x] \subseteq B \cup N_{G_{uv}}[x] \cup \{v\} = B \cup N_{G_{uv}}[x]$. Since the pair $(B,C)$ satisfies the \uvprop, so does $(B',C')$. Thus by induction, we note that
\[
\begin{array}{lcl}
1 + \dstartpass(G|C) & \stackrel{(x {\rm \, optimal \, move})}{=} & 1 + (1 + \sstartpass(G|C \cup N_G[x])) \1 \\
& \stackrel{({\rm by \, induction})}{\ge} & 1 + \sstart(G_{uv}|B \cup N_{G_{uv}}[x]) \2 \\
& \ge & \dstart(G_{uv}|B),
\end{array}
\]

\noindent
where the last inequality follows from the observation that the vertex $x$ may not be an optimal move for \D in the game played in $G_{uv}|B$. This completes Case~2.

\medskip
\emph{Case~3. $x \in \{u,v\}$ and $\xbar \notin B$.} Renaming $u$ and $v$ if necessary, we may assume that $x = u$ and $v \notin B$.  Suppose firstly that when \D plays the (optimal) vertex $x = u$, the only new vertex dominated in $G|C$ is the vertex $v$. As before, the pair $(B',C')$ satisfies the \uvprop. The following now holds.

\[
\begin{array}{lcl}
1 + \dstartpass(G|C) & \stackrel{(x \textrm{ optimal move})}{=} & 1 + (1 + \sstartpass(G|C \cup N_G[u])) \2 \\
& = & 1 + (1 + \sstartpass(G|C \cup \{v\}) \1 \\
& \stackrel{(\textrm{by induction})}{\ge} & 1 + \sstart(G_{uv}|B \cup \{v\}) \2 \\
& = & 1 + \sstart(G_{uv}|B \cup N_{G_{uv}}[u']) \2 \\
& \ge & \dstart(G_{uv}|B),
\end{array}
\]
\noindent
where the last inequality follows from the observation that the vertex $u'$ may not be an optimal move for \D in the game played in $G_{uv}|B$. Suppose secondly that when \D plays the (optimal) vertex $x = u$, at least one new vertex different from $v$ is dominated in $G|C$. In this case, $C \cup N_G[u] \setminus B \cup N_{G_{uv}}[u] = \{v\}$ and $\{u,v\} \subseteq C \cup N_G[u]$. As before, the pair $(B',C')$ satisfies the \uvprop\ and the following now holds, noting that $v \notin  N_{G_{uv}}[u]$.
\[
\begin{array}{lcl}
1 + \dstartpass(G|C) & \stackrel{(x {\rm \, optimal \, move})}{=} & 1 + (1 + \sstartpass(G|C \cup N_G[u])) \2 \\
& \stackrel{({\rm by \, induction})}{\ge} & 1 + \sstart(G_{uv}|B \cup N_{G_{uv}}[u]) \2 \\
& \ge & \dstart(G_{uv}|B),
\end{array}
\]
\noindent
where the last inequality follows from the observation that the vertex $u$ may not be an optimal move for \D in the game played in $G_{uv}|B$. This completes Case~1.
In all three cases above, we have shown that $\dstart(G_{uv}|B) \le 1 + \dstartpass(G|C)$. This completes the proof of Claim~B.
\smallqed

\medskip
Claim~A and Claim~B complete the proof of Lemma~\ref{lem:Cutting Lemma III}.
\qed

As a consequence of Theorem~\ref{thm:cutting-lemma}, Lemma~\ref{lem:Cutting Lemma II}, and Lemma~\ref{lem:Cutting Lemma III}, we now have the main result of this section.

\begin{thm}
{\rm\bf (Extended Cutting Lemma)}
\label{lem:cutting}
Let $G$ be a graph, and let $uv$ be an edge of $G$. If subsets of vertices $A$, $B$, and $C$ satisfy $C\subseteq B \subseteq A$, then the following hold.
\[
\begin{array}{rrcccl}
{\rm (a)} & \dstart(G|A) & \le & \dstart(G_{uv}|B) & \le & \dstartpass(G|C)+1. \2\\
{\rm (b)} & \sstart(G|A) & \le & \sstart(G_{uv}|B) & \le & \sstartpass(G|C)+1.
\end{array}
\]
\end{thm}

We remark that in order to prove the first inequality in both statements (a) and (b) of the Extended Cutting Lemma, we needed the stronger result involving the Staller-pass game in the three versions of the Cutting Lemma which allowed  us to prove the desired result on the domination game using induction. However we remark that when proving our main result in Section~\ref{sec:spiders}, we use the domination game result given by the Cutting Lemma and do not need the Staller-pass game result. From the Extended Cutting Lemma and Proposition~\ref{lem:pass-in-general}, we get the following consequence.

\begin{cor}
Let $G$ be a graph, and let $A,B,C \subseteq V(G)$ where $C \subseteq B\subseteq A$. If $uv$ is an edge of $G$, then the following hold.
\[
\begin{array}{rrcccl}
{\rm (a)} & \dstart(G|A) & \le & \dstart(G_{uv}|B) & \le & \dstart(G|C) + 2. \2\\
{\rm (b)} & \sstart(G|A) & \le & \sstart(G_{uv}|B) & \le & \sstart(G|C) + 2.
\end{array}
\]
\end{cor}

\section*{Acknowledgements}

Research of M.A.H. supported in part by the South African National Research Foundation and the University of Johannesburg. S.K.\, and G.K.\, acknowledges the financial support from the Slovenian Research Agency (research core funding No. P1-0297) and that the project (Combinatorial Problems with an Emphasis on Games, N1-0043) was financially supported by the Slovenian Research Agency. G.K.\ was also supported by Slovenian Research Agency under the grant P1-0294 and Slovenian Public Agency for Entrepreneurship, Internationalization, Foreign Investments and Technology under the grant KKIPP-99/2017.

\end{document}